\documentclass[11pt]{article}

\usepackage[english]{babel}
\usepackage{amsmath}
\usepackage{amssymb,amsfonts,a4}
\usepackage{eucal,bbm,stmaryrd}
\usepackage{psfrag,graphicx}

\begin{document}

\numberwithin{equation}{section}

\newtheorem{satz}{Theorem}
\newtheorem{defi}{Definition}
\newtheorem{lem}{Lemma}


\newcommand{\Om}{\Omega}
\newcommand{\Bew}{\noindent{\bf Proof: }}
\newcommand{\qed}{\hspace*{3em} \hfill{$\square$}}
\newcommand{\weg}{\setminus}
\newcommand{\gdw}{\Leftrightarrow}
\newcommand{\sgdw}{\Leftrightarrow}
\newcommand{\nach}{\longrightarrow}
\newcommand{\alle}{\,\forall\,}
\newcommand{\gibt}{\,\exists\,}
\newcommand{\Bed}{\,|\,}

\newcommand{\Z}{\mathbbm{Z}}
\newcommand{\R}{\mathbbm{R}}
\newcommand{\Q}{\mathbbm{Q}}
\newcommand{\baR}{\overline{\R}}
\newcommand{\RS}{\R^2_S}
\newcommand{\N}{\mathbbm{N}}

\newcommand{\p}{\mathcal{P}}
\newcommand{\prob}{\mathcal{P}_1}

\newcommand{\opB}{B_2}
\newcommand{\M}{\mathcal{M}}


\newcommand{\Y}{\mathcal{Y}}
\newcommand{\YY}{\mathcal{Y}_0}
\newcommand{\bY}{\bar{Y}}
\newcommand{\La}{\Lambda}
\newcommand{\Lan}{\La_n}

\newcommand{\Bo}{\mathcal{B}}
\newcommand{\B}{\mathcal{B}}
\newcommand{\F}{\mathcal{F}}

\newcommand{\la}{\lambda}
\newcommand{\si}{\sigma}

\newcommand{\einh}{\vec{e}}
\newcommand{\lear}{\le_{\einh}}
\newcommand{\Ga}{\Gamma}
\newcommand{\Gar}{\Ga_r}
\newcommand{\Gas}{\Ga_s}

\newcommand{\bs}{\bar{s}}

\newcommand{\E}{\mathcal{E}}

\newcommand{\gmu}{\gamma}
\newcommand{\G}{\mathcal{G}}


\newcommand{\gglatt}{\varphi}

\newcommand{\en}{E_{n}}

\newcommand{\ah}{\vec{a}} 
\newcommand{\tauh}{\vec{\tau}} 

\newcommand{\Nk}{\N_k}

\newcommand{\A}{\mathcal{A}}

\newcommand{\Gn}{G_{R,n}}
\newcommand{\Gnb}{G_n}
\newcommand{\Drn}{D_{R,n}}

\newcommand{\tn}{\tau_{R,n}}
\newcommand{\Taun}{\tau^q_{R,n}}
\newcommand{\bn}{\bar{n}}
\newcommand{\bR}{\bar{R}}
\newcommand{\bU}{\bar{U}}
\newcommand{\tU}{\tilde{U}}
\newcommand{\tiu}{\tilde{u}}

\newcommand{\pin}{\pi_n}
\newcommand{\pinh}{\pi'_n}

\newcommand{\myt}{m_{y',t}}
\newcommand{\hyt}{h_{y',t}}

\newcommand{\bpsi}{\bar{\psi}}

\newcommand{\de}{\delta}
\newcommand{\ep}{\epsilon}
\newcommand{\ex}{c_\xi}

\newcommand{\ph}{\varphi_{R,n}}
\newcommand{\bph}{\bar{\varphi}_{R,n}}
\newcommand{\leb}{\lambda^2}
\newcommand{\cK}{c_K}
\newcommand{\cpsi}{c_{\psi}}
\newcommand{\cu}{c_u}
\newcommand{\KU}{K^{U}}
\newcommand{\Ko}{K}
\newcommand{\Kop}{K'}
\newcommand{\Kopp}{K''}
\newcommand{\Kepa}{K_{\ep,\ah}}
\newcommand{\fc}{f_{\Ko}}
\newcommand{\cf}{c_f}
\newcommand{\ft}{\tilde{f}}

\newcommand{\per}{\eta}
\newcommand{\akb}{A_{k,B,\per}}
\newcommand{\takb}{\tilde{A}_{k,B,\per}}

\newcommand{\ela}{\Lambda}


\newcommand{\ayb}{{\tn^{\wedge,Y,B_+}}}
\newcommand{\ryb}{r^{Y,B_+}_{n'}}

\newcommand{\T}{\mathfrak{T}}
\newcommand{\Tn}{\T_{R,n}}
\newcommand{\Tny}{\T_{R,n}^{Y}}
\newcommand{\Tnb}{\T_{R,n}^{B}}

\newcommand{\tyb}{t^{Y,B}}
\newcommand{\Tyb}{T^{Y,B}}
\newcommand{\tauyb}{\tau^{Y,B}}
\newcommand{\pyb}{P^{Y,B}}
\newcommand{\cyb}{C^{Y,B}}

\newcommand{\etyb}{t}
\newcommand{\eTyb}{T}
\newcommand{\etauyb}{\tau}
\newcommand{\epyb}{P}
\newcommand{\ecyb}{C}

\newcommand{\tY}{\tilde{Y}}
\newcommand{\tB}{\tilde{B}}
\newcommand{\ty}{\tilde{y}}
\newcommand{\tm}{\tilde{m}}
\newcommand{\tP}{\tilde{P}}
\newcommand{\tC}{\tilde{C}}
\newcommand{\ttau}{\tilde{\tau}}
\newcommand{\tT}{\tilde{T}}

\newcommand{\tTn}{\tilde{\T}_{R,n}}
\newcommand{\tTny}{\tilde{\T}_{R,n}^{\tY}}
\newcommand{\tTnb}{\tilde{\T}_{R,n}^{\tB}}

\newcommand{\ttyb}{\tilde{t}^{\tY,\tB}}
\newcommand{\tTyb}{\tilde{T}^{\tY,\tB}}
\newcommand{\ttauyb}{\tilde{\tau}^{\tY,\tB}}
\newcommand{\tpyb}{\tilde{P}^{\tY,\tB}}
\newcommand{\tcyb}{\tilde{C}^{\tY,\tB}}

\newcommand{\ettyb}{\tilde{t}}
\newcommand{\etTyb}{\tilde{T}}
\newcommand{\ettauyb}{\tilde{\tau}}
\newcommand{\etpyb}{\tilde{P}}
\newcommand{\etcyb}{\tilde{C}}

\newcommand{\iTn}{\bar{\T}_{R,n}}
\newcommand{\iTny}{\bar{\T}_{R,n}^{Y}}
\newcommand{\iTnB}{\bar{\T}_{R,n}^{B}}
\newcommand{\iTyb}{\bar{T}_{R,n}^{Y,B}}
\newcommand{\itny}{\bar{t}_{R,n}^{Y,B}}

\newcommand{\tla}{\tilde{\la}}


\begin{center}
{\bf \Large Translation-invariance of two-dimensional
Gibbsian systems of particles with internal degrees of freedom}\\

\vspace{.75 cm}

Thomas Richthammer\\
Department of Mathematics, UCLA, Los Angeles, CA 90095-1555\\
Email: richthammer@math.ucla.edu \\
Tel: +1 310 206 8870, Fax: +1 310 206 6673
\end{center}
\vspace{ -0.1 cm}

\begin{abstract}
One of the main objectives of equilibrium state statistical physics is to analyze which symmetries of an interacting particle system in equilibrium
are broken or conserved. Here we present a general result on the conservation of translational symmetry for two-dimensional Gibbsian particle systems. The result applies to particles with internal degrees of freedom and fairly arbitrary interaction, including the interesting cases of discontinuous, singular, and hard core interaction.  In particular we thus show the conservation of translational symmetry for the continuum Widom Row\-linson model and a class of continuum Potts type models. 
\\
 
\end{abstract}


\section{Introduction}

\begin{sloppypar}

It is well known that probability theory provides a mathematically rigorous setting to investigate problems from equilibrium state statistical physics. Here the object of consideration is a system of interacting particles, where the number of particles is huge and thus assumed to be infinite. Such a particle system is given by specifying  restrictions on particle positions (lattice setting versus point particle setting), the internal properties of the particles
(such as magnetic spin, electric charge or particle type), and 
the interaction between particles.  The equilibrium states of a specific particle system are then modeled by Gibbsian processes, as introduced by  R.~L.~Dobrushin (see \cite{D68} and \cite{D70}), O.~E.~Lanford and D.~Ruelle (see \cite{LR}). 
The main objective usually is to find out whether the system exhibits a phase transition, i.e. whether there is more than one equilibrium state. In order to study this problem, the crucial task is to investigate which of the system's symmetries are broken and which are conserved, and it would be desirable to have general results stating under which conditions certain symmetries are conserved. While in more than two spatial dimensions such general results can not be expected to hold (as here all symmetries are believed to be broken easily), and in one dimension the situation is almost trivial (as under very weak assumptions all symmetries are conserved), the case of two dimensions is interesting. 
Here it is useful to distinguish between discrete and continuous symmetries and also between internal symmetries (i.e. symmetry transformations concerning the inner properties of particles) and spatial symmetries (such as translation and rotation). In order to investigate the behavior of a particle system under translations and rotations it is natural to consider a point particle setting,
as in a lattice setting all spatial symmetries are bound to be discrete. In the following the attention is thus restricted to  interacting particle systems in a point particle setting in two dimensions. \\

The current knowledge about such systems is the following: It is believed that discrete internal symmetries in general may be broken, but so far this has been shown only for very few systems, e.g. the Widom Rowlinson model considered by D.~Ruelle \cite{Ru2} or the continuum Potts model considered by H.-O.~Georgii and O.~H\"aggstr\"om \cite{GH}. In contrast, continuous internal symmetries are conserved under weak assumptions on the interaction. The first result in this direction was obtained by S.~Shlosman \cite{S}, building on earlier ideas of M.~Mermin and H.~Wagner \cite{MW}. We gave a more general version of this result in \cite{Ri1}, which includes the case of discontinuous interaction, using ideas of  D.~Ioffe, S.~Shlosman and Y.~Velenik \cite{ISV}. While it is expected that rotational symmetry may be broken, so far this could not be established for any realistic particle system, but this conjecture is supported by recent work of F.~Merkl and S.~Rolles \cite{MR}, for example. Translational symmetry is conserved under weak assumptions on the interaction. This was first shown by E.~Fr\"ohlich and C.-E.~Pfister \cite{FP1} and \cite{FP2}, and we obtained a more general result \cite{Ri2}, which for example includes the interesting case of the hard disc model. The last two results both concern particles without internal degrees of freedom.  However, many interesting models of statistical physics, such as the Widom Rowlinson or the Potts model, feature particles with spins. Here we will show, how to overcome conceptual and technical difficulties that arise due to the incorporation of spins, and thus we obtain a fairly general result on the conservation of translational symmetry for particles with any internal degrees of freedom and for interactions that are allowed to have discontinuities, singularities, or hard cores. This establishes the conservation of translational symmetry for the continuum Widom Rowlinson model and a large class of continuum Potts type models, for example. While parts of the proof of
the main theorem will be similar to the corresponding parts in \cite{Ri2}, we decided to repeat these arguments for the convenience of the reader, so that the article is self-contained.\\ 

We start Section~\ref{secresult} by giving an equivalent condition 
for a measure to be invariant under a transformation (Lemma~\ref{lekrit}), which will be useful for establishing the conservation of symmetries. We next confine ourselves to the special case of translational symmetry, considering a class of 
Potts type potentials. The corresponding result 
(Theorem~\ref{pottstype}), which is of interest on its own, will follow from 
the general case presented afterwards. For this general case we define a class 
of potentials (Definition~\ref{defGapprox}) for that translational symmetry is 
conserved  (Theorem~\ref{sym}). After a few comments on some aspects of 
this class concerning hard cores (Lemmas~\ref{lebedeperw} and \ref{leballs}) we  give sufficient conditions for potentials to belong to this class 
(Lemmas~\ref{lepotglatt} and \ref{lepotpotts}). The precise setting is given in Section~\ref{secsetting}, and 
the proofs of the lemmas from Sections~\ref{secresult} and \ref{secsetting} are
relegated to Section~\ref{secleset}. In Section~\ref{secproofsym} we will 
give the proof of Theorem~\ref{sym}. The proofs of the corresponding lemmas are 
relegated to Section~\ref{secleproofsym}.

%

\section{Results} \label{secresult}

\subsection{Conservation of symmetries}

We consider particles in the plane $\R^2$. Every particle is allowed to have  
internal degrees of freedom, encoded in the so called spin of the particle. 
The spin is assumed to be an element of some measurable spin space 
(or mark space) $(S,\F_S)$, on which a probability measure $\la_S$ is given as a 
reference measure. We require the diagonal in $S \times S$
to be measurable w.r.t $\F_S \otimes \F_S$, but we will not assume any topological properties of $S$. The particle space will be abbreviated 
by $\RS := \R^2 \times S$. 
We fix a chemical potential $-\log z$, where $z>0$ is a given 
activity parameter. The particles may interact via a \emph{pair potential} 
$U$ modelled by a measurable function 
\[
U: (\RS)^2 \; \to \; \baR \; := \; \R \cup \{\infty\}
\]
that is symmetric in that $U(y_1,y_2) = U(y_2,y_1)$. The set of 
equilibrium states corresponding to a particular choice of $U$ and $z$ 
can be modelled by the set of Gibbsian point processes, which are 
defined to be certain probability measures on the space  $(\Y,\F_{\Y})$ 
of all particle configurations, see Section~\ref{secgibbs}.
A bimeasurable transformation $\tau: \RS \to \RS$ 
is called a \emph{symmetry} of $U$ if $U$ is $\tau$-invariant, i.e. 
\begin{equation*}
U(\tau (y_1), \tau (y_2)) = U(y_1,y_2) \quad 
\text{ for all } y_1,y_2 \in \RS.
\end{equation*}
Such a transformation $\tau$ also defines a transformation on the configuration space $\Y$,  where every single particle of a given configuration is transformed by $\tau$, and an equilibrium state $\mu$ is said to be $\tau$-invariant if $\mu \circ \tau^{-1} = \mu$. 
It is natural to ask, whether the equilibrium states 
of a particle system corresponding to $U$ and $z$ are invariant 
under a given symmetry of $U$. If this is indeed the case, 
the symmetry is said to be conserved, otherwise it is said to be broken.\\

There are several strategies to establish the conservation of symmetries. 
One is to use the concept of relative entropy and to exploit 
certain entropy estimates, see Section~2.3.3. of \cite{ISV}. 
Another one builds on a certain inequality for Gibbsian specifications, 
see Proposition (9.1) of \cite{G}. The latter approach uses the 
convexity of the set of Gibbs measures, tail triviality of extremal Gibbs 
measures and extreme decomposition, thus requiring the spin space to be 
standard Borel. In the following we present a variant of 
this approach, which works in a general setting and admits a straightforward 
proof via convexity. 

\begin{lem} \label{lekrit}
Let $(\Om,\F,\mu)$ be a probability space, $\A \subset \F$ an algebra 
on $\Om$ such that $\si(\A)=\F$, and $\tau$ a transformation
on $\Om$, i.e.  $\tau: \Om \to \Om$ a  bimeasurable mapping. 
$\mu$ is $\tau$-invariant if and only if the following condition holds: 
\begin{equation} \label{kritA} 
\forall A \in \A: \quad 
\mu(\tau A) \, + \, \mu(\tau^{-1} A) \, \ge \, 2 \mu(A). 
\end{equation}
\end{lem}
The proof will be given in Section~\ref{secleset}. For a more detailed 
account on how to use this lemma in order to show the conservation 
of translational symmetry, see Subsection~\ref{secsym}. From now on we 
will restrict our attention to spatial translations of particles. 


\subsection{Widom Rowlinson and Potts type potentials} \label{secpotts}

As the definition of the class of potentials for that we will show 
the conservation of translational symmetry is fairly general, but 
also fairly complicated, we first would like to present the result 
for a certain class of Potts type potentials. This class includes 
finite state Widom Rowlinson potentials as wells as step potentials, 
as considered by J.~L.~Lebowitz and E.~H.~Lieb in \cite{LL} as a continuum 
analogue of the Potts model. For a given finite spin space $S$ 
(describing different types of particles) we define a \emph{Potts type 
potential} to be of the form 
\[
U(x_1,\si_1,x_2,\si_2) := \phi_{\si_1,\si_2}(|x_1-x_2|_h),
\]
where $|.|_h$ is a norm on $\R^2$ and $(\phi_{\si_1 \si_2})_{\si_1,\si_2 \in S}$ 
is a family of interactions, i.e. $\phi_{\si_1 \si_2}: \R_+ := [0,\infty[ \to \baR$ is  
measurable and we have $\phi_{\si_1 \si_2} = \phi_{\si_2 \si_1}$ 
for all $\si_1,\si_2$. 
\begin{figure}[!htb] 
\begin{center}
\psfrag{pa}{$\phi_a(r)$}
\psfrag{pb}{$\phi_b(r)$}
\psfrag{pc}{$\phi_c(r)$}
\psfrag{r}{$r$}
\psfrag{i}{$\infty$}
\psfrag{1}{$1$}
\includegraphics[scale=0.4]{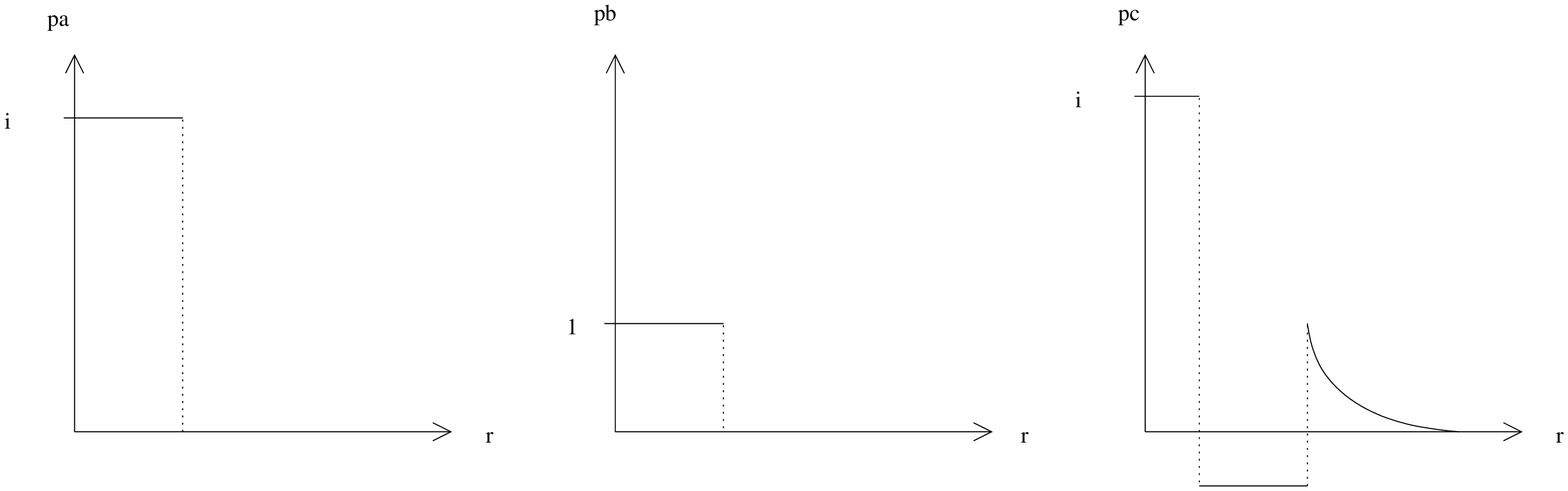}
\end{center}
\caption{Some examples of well behaved functions} \label{figpot}
\end{figure}

\noindent
We call a function $\phi:\R_+ \to \baR$ \emph{well behaved} if 
there are $0 \le r_0 < \ldots < r_n$  $(n \ge 0)$ such that 
$\phi(r) = \infty$ for $r < r_0$, $\phi(r) = 0$ for $r > r_n$,  
$\phi$ is continuous on every interval $]r_i,r_{i+1}[$ and in 
every point $r_1, \ldots, r_n$ the left and right limit exist. 
Figure \ref{figpot} shows some examples of well behaved functions. Functions of 
type $\phi_a$, $\phi_b$ and $\phi_c$ are used in the definition of a Widom Rowlinson 
potential, a continuum Potts potential and a slightly more complicated 
Potts type potential respectively.
All spatial translations of particles are symmetries of Potts type potentials,
and the following theorem states the conservation of these symmetries. 
\begin{satz} \label{pottstype}
Let $S$ be a finite spin space endowed with the equidistribution as 
reference measure and $z>0$ be an activity  parameter. Let 
$U$ be a Potts type potential corresponding to a norm $|.|_h$ on $\R^2$ 
and a family of interactions $(\phi_{\si_1 \si_2})_{\si_1,\si_2 \in S}$. 
If all the functions $\phi_{\si_1 \si_2}$ are nonnegative and well behaved, 
then every Gibbs measure corresponding to $U$ and $z$ is translation-invariant.
\end{satz}   
We note that nonnegativity of the potential is assumed only in order to 
avoid introducing superstability at this point. In Section~\ref{secsmap}, Theorem~\ref{pottstype} will be deduced from the general case (Theorem~\ref{sym}) presented below. 


\subsection{General case}

In this general case we consider translations in a fixed direction
$\ah$ ($\ah \in \R^2$ with $|\ah|_2 = 1$). 
The corresponding group of translation transformations is defined by 
\[
g_{t}: \RS \to \RS, \quad 
g_{t}(x,\si) := (x,\si) + \ah t := (x+\ah t,\si) \quad (t \in \R). 
\]
We call a potential $U$ (or a  Gibbsian point process $\mu$) 
invariant under translations in direction $\ah$ or simply 
$\ah$-invariant if $U$ (or $\mu$ respectively) is invariant 
under $g_t$ for all $t \in \R$. 
Translation-invariance is equivalent to $\ah$-invariance in 
every direction $\ah$. As there might be interesting potentials
that are $\ah$-invariant for some direction $\ah$, but not
for every direction, we investigate the conservation of
$\ah$-translational symmetry rather than translational symmetry.\\ 

In order to describe a class of potentials for that $\ah$-symmetry 
is conserved, we now define some important properties of sets, functions, 
and potentials. We call a function $f:  (\RS)^2 \to \baR$ 
\[
\begin{array}{ll}
\text{\emph{$\ah$-invariant} }  
&\text{ if } \; f(y + t\ah,y' + t\ah) = f(y,y') \;  
 \alle y,y' \in \RS, t \in \R,\\
\text{\emph{symmetric} } 
&\text{ if } \; f(y,y') =f(y',y) \;  \alle y,y' \in \RS \quad  \text{ and }\\
\text{of \emph{bounded range}}
&\text{ if } \; \{ |y-y'|: f(y,y') \neq 0\} \; \text{ is bounded.} 
\end{array}
\]
Here the distance of two particles is defined to be the distance
of the positions of the particles. The above definition of course does not
depend on the choice of norm $|.|$, but for sake of definiteness let
$|.|$ be the maximum norm on $\R^2$. We say that a set $A \subset  (\RS)^2$
is $\ah$-invariant, symmetric, or of bounded range if the corresponding indicator 
function $1_A$ has this property. We call $A$ a \emph{standard set}
if it is measurable, symmetric, and of bounded range. 
Let us call $U$ a \emph{standard potential} if it is measurable, symmetric, 
and its  hard core 
\[
\KU \, := \, \{ U = +\infty\}
\] 
is a standard set, i.e. if its hard core is of bounded range. 
Usually the hard core can be described in terms of a norm,
which is the case for Potts type potentials as described above, 
for example, but in our setup we are able to treat 
fairly general hard cores.
We also need regularity properties. We call a potential $U$ 
\emph{$\ah$-continuous}, \emph{$\ah$-equicontinuous}, or \emph{$\ah$-smooth} 
on a set $A$ if the family of functions
\[
\gglatt^U_{y_1,y_2}(t): \R \to \baR, \quad t \mapsto U(y_1,y_2 + t\ah) 
\quad ((y_1,y_2) \in A)
\]
is  continuous, equicontinuous, or smooth in $t = 0$. In the case of 
smoothness we define the $\ah$-derivatives of $U$ in $A$ by 
\[
\partial^2_{\ah} U(y_1,y_2) \, := \, \frac{d^2}{dt^2} \gglatt^U_{y_1,y_2}(0), 
\] 
and for a given function  $\psi:(\RS)^2 \to \R_+$ 
we say that the $\ah$-derivatives of $U$ are \emph{dominated} by $\psi$ on $A$ if   
\[
\partial^2_{\ah} U(y_1,y_2+ t\ah) \le \psi(y_1,y_2)
\; \text{ for all } (y_1,y_2) \in A, t \in [-1,1] \text{ s.t. }  
(y_1,y_2+ t\ah) \in A.
\]
In the context of $\psi$-domination we will use the notion 
of a \emph{bounded partially square integrable function (bpsi-function)}, 
which is defined to be a measurable, symmetric function  $\psi:(\RS)^2 \to \R_+$  
satisfying
\[
\|\psi\| \, < \, \infty  \quad \text{ and } \quad 
\sup_{y_1 \in \RS} \int \psi(y_1 ,y_2) |y_1 -y_2|^2  dy_2 \, < \, \infty,
\]
where $\|.\|$ is the supremum norm of a function. In order to be able to control
the potential in a neighborhood of a given set, we introduce the notion of
the \emph{$\ep$-$\ah$-enlargement} $K_{\ep,\ah}$ of a set  
$K \subset (\RS)^2$ for a given $\ep>0$, defined by 
\[
\Kepa  \, := \,  \{(y_1,y_2 + r \ah): (y_1,y_2) \in K, - \ep < r < \ep\}. 
\]
We note that the $\ep$-$\ah$-enlargement of an $\ah$-invariant standard set 
again is an $\ah$-invariant, symmetric set of bounded range. However,
it is not necessarily measurable, so we need to be a bit more careful. 
Given an $\ah$-invariant standard set $\Ko$ 
we define $\Kop$, $\Kopp$ to be \emph{measurable $\ah$-enlargements of $\Ko$} if for some $\ep>0$ 
\[
\text{$\Kop$ and $\Kopp$ are standard sets, } \quad 
\Ko_{\ep,\ah} \subset \Kop \quad \text{ and } \quad (\Kop)_{\ep,\ah}\subset \Kopp.
\]
If $U$ is a potential, $z>0$ is an activity parameter and $\YY$ is a
set of boundary conditions, we say that the triple $(U,z,\YY)$ is
admissible if all conditional Gibbs distributions corresponding to $U$ and $z$ 
with boundary condition taken from $\YY$ are well defined, see
Definition~\ref{defzul} in Section~\ref{secgibbs}. Important examples
are the cases of superstable potentials with tempered boundary configurations 
and nonnegative potentials with arbitrary boundary conditions, 
see Section~\ref{secadm}. For admissible  $(U,z,\YY)$ the set of Gibbs measures 
$\G_{\YY}(U,z)$ corresponding to $U$ and $z$ with full weight on configurations 
in $\YY$ is a well defined object. Finally we need bounded correlations: 
For admissible $(U,z,\YY)$ we call $\xi \in \R$ a Ruelle bound if the 
correlation function of every Gibbs measure $\mu \in \G_{\YY}(U,z)$ is bounded 
by powers of $\xi$ in the sense of \eqref{ruellebd} in Section~\ref{secgibbs}.
\begin{defi} \label{defGapprox}
Let $(U,z,\YY)$ be admissible with Ruelle bound $\xi$, 
where $U:  (\RS)^2 \to \baR$ is an $\ah$-invariant 
standard potential. We say that $U$ is $\ah$-smoothly approximable
if there is a decomposition of $U$ into a smooth part $\bU$ and
a small part $u$ in the following sense:  
We have an $\ah$-invariant standard set $\Ko \supset \KU$ and 
measurable symmetric $\ah$-invariant functions $\bU, u: \Ko^c \to \R$
such that $U =  \bU - u$, $u  \ge  0$, $\bU$ has  $\psi$-dominated $\ah$-derivatives on  $\Ko^c$ for some bpsi-function $\psi$, and 
\begin{equation} \label{Gapprox}
\begin{split}
&\sup_{y_1 \in \RS} \int \,(1_{\Ko^c} \tiu) (y_1,y_2)  
|y_1-y_2|^2 \, dy_2 \,  < \, \infty  \, \text { and } \\
&\sup_{y_1 \in \RS}  \int  \, (1_{\Kopp \weg \KU} + 1_{\Ko^c}\tiu) 
(y_1,y_2)\,dy_2 \, < \, \frac 1 {z\xi} 
\end{split}
\end{equation}
for some  measurable $\ah$-enlargements $\Kop$, $\Kopp$ of $\Ko$ and 
$\tiu := 1 - e^{-u} \le  u \wedge 1$.
\end{defi}
\pagebreak
The class of smoothly approximable standard potentials is a rich 
class of potentials. An $\ah$-smoothly approximable $\ah$-invariant 
standard potential may  have a singularity or a hard core at the origin, 
and the type of convergence into the singularity or the hard core is 
fairly arbitrary, as we have not  imposed any condition on $U$ in 
$\Ko \weg \KU$. For small activity $z$ the last condition of \eqref{Gapprox} 
holds for large sets $\Kopp$, so $\Ko$ can be chosen to be a large set, 
which relaxes the conditions on $U$. The small part $u$ of $U$ is not 
assumed to satisfy any regularity conditions, so that $U$ doesn't have to be 
smooth or continuous. We note that Definition~\ref{defGapprox} 
does not depend on the choice of the norm $|.|$.\\

The above definition may seem overly complicated. Nevertheless we 
present it in the given form in order to include as many potentials 
as possible in the class of $\ah$-smoothly approximable potentials. 
For some comments on Definition~\ref{defGapprox} and simplifications 
in several special cases we refer to the following subsections. Beforehand 
however, we would like to present our main result: 
\begin{satz} \label{sym}
Let $(U,z,\YY)$ be admissible with Ruelle bound, 
where $U: (\RS)^2 \to \baR$ is a standard potential. 
If $U$ is $\ah$-invariant and $\ah$-smoothly approximable, then 
every Gibbs measure $\mu \in \G_{\YY}(U,z)$ is $\ah$-invariant. 
\end{satz}


\subsection{Measurable enlargements} \label{secmeasen}

In Definition~\ref{defGapprox} the hard core $\KU$ may be of fairly arbitrary 
size and shape. The only condition on $\KU$ is the existence of an 
$\ah$-invariant standard set $\Ko \supset \KU$ and of measurable $\ah$-enlargements  $\Kop$, $\Kopp$ of $\Ko$ such that  $\Kopp \weg \Ko$ and $\Ko \setminus \KU$ are not too big in the sense of the second inequality of \eqref{Gapprox}. In the following we present 
possible constructions of  measurable $\ah$-enlargements $\Kop$, $\Kopp$ 
for a given $\ah$-invariant standard set $\Ko$.\\

{\bf (a)} Even if the $\ep$-$\ah$-enlargements of a measurable  set $\Ko$ in general need not be measurable again, they often are. If so, 
we simply set $\Kop = \Ko_{\ep,\ah}$ and  $\Kopp=(\Kop)_{\ep,\ah}=\Ko_{2\ep,\ah}$
to construct measurable enlargements of $\Ko$.
\begin{lem} \label{lebedeperw}
Let $A \subset (\RS)^2$ be a measurable set with 
all $\ah$-cross sections
$~A(y_1,y_2,\ah) := \{r \in \R: (y_1,y_2+r\ah) \in A\}$ $~(y_1,y_2 \in \RS)$  satisfying
\begin{equation} \label{bedeperw}
\alle U \subset \R \text{ open}: \quad 
U \cap A(y_1,y_2,\ah) \neq \emptyset \,  \Rightarrow \, 
\la^2(U \cap A(y_1,y_2,\ah)) > \, 0. 
\end{equation}
Then every  $\ep$-$\ah$-enlargement $A_{\ep,\ah}$ $\,(\ep > 0)$ is measurable again. 
For example \eqref{bedeperw} is satisfied if the set of interior points 
of $A(y_1,y_2,\ah)$ is dense in  $A(y_1,y_2,\ah)$.   
\end{lem}   
Condition \eqref{bedeperw} concerns only the topological structure 
of the $\ah$-cross sections of $A$,  is easy to be verified, and 
holds in the case that each $\ah$-cross section of $A$ 
is open, for example. \\

{\bf (b)}
Often we may choose $\Ko$ to \emph{consist of discs} 
in the following sense: 
\[
\Ko = \{ (x_1,\si_1,x_2,\si_2) \in (\RS)^2 : 
|x_1-x_2|_h \le r_{\si_1 \si_2}\}, 
\]
where $|.|_h$ is an arbitrary norm on $\R^2$ and 
$(r_{\si_1 \si_2})_{\si_1,\si_2 \in S}$ is a symmetric, measurable and bounded family of reals. In this case we define the enlargement 
$\Ko_{+\ep}$ to be a set of the above form, 
where $r_{\si_1 \si_2}$ is replaced by $r_{\si_1 \si_2} + \ep$. 
We now simply set $\Kop := \Ko_{+\ep}$ and $\Kopp := \Kop_{+\ep} = \Ko_{+2\ep}$.\\

{\bf (c)} 
If $(S,\F_S)$ is assumed to be a standard Borel space, then
there is a metric on $S$, such that $\F_S$ is the corresponding 
Borel-$\si$-Algebra. Hence there is a measurable metric $d$ on 
$(\RS)^2$, and the $d,\ep$-enlargement of an arbitrary set 
$A \subset (\RS)^2$ 
\[
A_{d,\ep}\, := \,\{ (y_1,y_2) \in  (\RS)^2: d((y_1,y_2),A) < \ep\}
\]
is measurable. In this case we may set $\Kop = \Ko_{d,\ep}$ and  $\Kopp=(\Kop)_{d,\ep}=\Ko_{d,2\ep}$. \\

We note that in the above cases we often may replace $\Kopp$ by $\Ko$ in the second inequality of \eqref{Gapprox}. Here the following easy lemma
is useful:
\begin{lem} \label{leballs}
Let $\Ko \subset (\RS)^2$ be an $\ah$-invariant standard set. 
For $\ep \to 0$
\begin{enumerate}
\item[(a)] 
$\sup_{y_1 \in \RS} 
\int 1_{\Ko_{\ep,\ah} \weg \Ko}(y_1,y_2) dy_2 \to 0$ 
in the situation of case (a) if we know that all $\ah$-cross
section of $\Ko$ are open intervals,
\item[(b)]
$\sup_{y_1 \in \RS} 
\int 1_{\Ko_{+\ep} \weg \Ko}(y_1,y_2) dy_2 \to 0$ 
in the situation of case (b).
\end{enumerate}
\end{lem}


\subsection{Smoothly approximable potentials} \label{secsmap}

For convenience, in the following we will stick to the case
of a hard core consisting of discs.  
The results show, how Definition~\ref{defGapprox} 
simplifies whenever $U$ satisfies additional regularity properties
such as smoothness or continuity.
\begin{lem} \label{lepotglatt}
Let $(U,z,\YY)$ be admissible with Ruelle bound, 
where $U: (\RS)^2 \to \baR$ is an $\ah$-invariant 
standard potential with a hard core $\KU$ consisting of discs. 
$U$ is $\ah$-smoothly approximable if 
for every $\de > 0$ there is a bpsi-function $\psi_{\de}$ such that  
one of the following conditions is satisfied: 
\begin{enumerate}
\item[(a)]
$U$ is smooth and has $\psi_{\de}$-dominated $\ah$-derivatives 
on $(\KU_{+\de})^c$. 
\item[(b)]
$U$ is bounded and $\ah$-equicontinuous on $(\KU_{+\de})^c$, 
and for some $R \in \N$ the set $\tilde{K}  :=  \{(y_1,y_2) \in (\RS)^2: |y_1 -y_2| \le R\}$
has one of the following properties:
\begin{enumerate}
\item[(b1)] $U$ has $\psi_{\de}$-dominated $\ah$-derivatives on $\tilde{K}^c$. 
\item[(b2)] There is an $\ah$-invariant standard potential $\tU$ such that 
$|U| \le \tU$ on $\tilde{K}^c$, $\tU$ is bounded and has $\psi_{\de}$-dominated $\ah$-derivatives on $\tilde{K}^c$,  and 
\[
\lim_{r \to \infty} \sup_{y_1 \in \RS}  
 \int \tU(y_1,y_2) |y_1-y_2|^2 1_{\{|y_1-y_2| \ge r\}}dy_2 \, = \, 0.
\] 
\end{enumerate}
\end{enumerate}
\end{lem}
\pagebreak
For example, $(b1)$ holds trivially when $U$ has finite range, 
and $(b2)$ includes the case that there are  $\ep' > 0$ and $k \ge 0$ 
such that $|U(y_1,y_2)| \le k / |y_1-y_2|^{4+\ep'}$ for large $|y_1-y_2|$.
We note that the generalization of the preceding lemmas  to 
more general hard cores is straightforward: 
Instead of imposing regularity conditions for $U$ on $(\KU_{+\de})^c$ 
for every $\de>0$ we just require it on  some $\ah$-invariant 
standard set $\Ko \supset \KU$ such that there are measurable 
$\ah$-enlargements  $\Kop$, $\Kopp$ of $\Ko$ with 
$\int 1_{\Kopp \weg \KU}(y_1,y_2)dy_2 < 1/(z\xi)$ for all 
$y_1 \in \RS$, where $\xi$ is a Ruelle bound.
Finally we would like to show how to deal with discontinuous potentials 
by considering Potts type potentials as defined in Section~\ref{secpotts}. 
\begin{lem} \label{lepotpotts}
Let $S$ be a finite spin space and let $(U,z,\YY)$ be  admissible 
with Ruelle bound, 
where $U$ is a Potts type potential corresponding to a norm $|.|_h$ on $\R^2$ 
and a family of interactions $(\phi_{\si_1 \si_2})_{\si_1,\si_2 \in S}$. 
If all the functions $\phi_{\si_1 \si_2}$ are well behaved, then 
$U$ is $\ah$-smoothly approximable.
\end{lem}
Again we note that the ideas of the proof of Lemma~\ref{lepotpotts} can 
be used to prove $\ah$-smooth approximability for more general potentials, 
but for simplicity we restrict ourselves to the above case. Theorem~\ref{pottstype} can now be seen to be an immediate consequence of 
Theorem~\ref{sym} and Lemma~\ref{lepotpotts}. We only have to note 
that for nonnegative $U$  $~(U,z,\Y)$ is admissible and admits a Ruelle bound (see Section~\ref{secadm}, Lemma~\ref{lemsup}). 

%

\section{Setting} \label{secsetting}

\subsection{State space}

We will use the notations $\N := \{0,1,\ldots\}$, 
$\R_+ := [0,\infty[$, $\baR := \R \cup \{+\infty\}$,  
$r_1 \vee r_2  :=  \max\{r_1,r_2\}$, and 
$r_1 \wedge r_2  :=  \min\{r_1,r_2\}$  for $r_1,r_2 \in \R$.
For sets $A,B$ the cross sections of a subset $C \subset A \times B$ with 
respect to given elements $a \in A$ and $b \in B$ are denoted by  
\[
C(a) \; := \; \{ b' \in B: (a,b') \in C\} \quad \text{ and } \quad 
C(b) \; := \; \{ a' \in A: (a',b) \in C\}. 
\] 
The state space  $\RS := \R^2 \times S$ of a particle consists of the space
of positions $\R^2$ and the spin space $S$. Usually we will denote
particles by $y$, positions by $x$ and spins by $\si$. 
Considering a model that does not include internal properties of 
particles we may simply set $S := \{0\}$. On  $\R^2$ we consider the 
maximum norm $|.|$ and the Euclidean norm $|.|_2$. 
The  Borel-$\si$-algebra $\B^2$ on $\R^2$
is induced by any of these norms. Let $\B_b^2$ be the set of all 
bounded Borel sets and $\la^2$ be the Lebesgue measure on  $(\R^2,\Bo^2)$. 
Integration with respect to this measure will be abbreviated by 
$dx := d\la^2(x)$. Often we consider the centred squares 
\[
\La_r \,:= \,[-r,r[^2 \, \subset \, \R^2 \qquad (r \in \R_+).
\]
For describing the spins of the particles let $(S,\F_S,\la_S)$
be a probability space such that the diagonal in $S \times S$
is measurable w.r.t $\F_S \otimes \F_S$. 
Integration with respect to $\la_S$ will be 
abbreviated by $d\si := d\la_S(\si)$ and the same way we use 
$dy := d\leb(x)d\lambda_S(\sigma)$ in the particle space.
Sometimes we will apply functions of $\R^2$ 
to particles by simply ignoring their spins; for example 
$|y_1-y_2|$ is defined to be the distance between the positions of
two particles $y_1,y_2 \in \RS$. Similarly  we may think of
$\La \subset \R^2$ as a set of particles by identifying this set with 
$\La \times S \subset \RS$.\\

We also want to consider bonds between particles. 
For a set $X$ we define 
\[
E(X)  \, :=  \, \{ A \subset X: \#A = 2\}
\] 
to be the set of all bonds in $X$, 
where $\#$ denotes the cardinality of a set. 
A bond will be denoted by $xx' :=\{ x,x'\}$, 
where $x,x' \in X$ such that $x \neq x'$. Every symmetric function $u$ 
on $X \times X$ can be considered a function on $E(X)$ via 
$u(xx') := u(x,x')$. For a bond set $B \subset E(X)$  $~(X,B)$ 
is a (simple) graph. The connectedness relation 
\[
x \stackrel{X,B}{\longleftrightarrow} x' \; :\gdw \; 
\gibt m \in \N, x_0,\ldots,x_m \in X: \, 
x = x_0, x' = x_m, \, 
x_{i-1}x_{i} \in B \; \alle i 
\]
defines an equivalence relation on $X$ whose equivalence classes
are called the $B$-clusters of $X$. Let  
\[
C_{X,B}(x) \, :=\, \{ x' \in  X: x \stackrel{X,B}
                                {\longleftrightarrow} x' \}   
 \quad \text{ and } \quad 
 C_{X,B}(\La) \, := \,  \bigcup_{x' \in X \cap \La} C_{X,B}(x') 
\]
denote the $B$-clusters of a point $x$ and a set $\La$ respectively. 
Primarily we are interested in the case  $X \subset \RS$. 
On $E(\RS)$ we consider the 
$\si$-algebra \[
\F_{E(\RS)} \, := \, \{ \{ y_1y_2 \in E(\RS):  
(y_1,y_2) \in M \} :  \, M \in (\B^2 \otimes \F_S)^2 \}.
\]


\subsection{Configuration space}

A set of particles $Y \subset \RS$ is called  
\[
\begin{array} {ll}
\text{finite } \quad &\text{ if } \#Y < \infty, \qquad \text{ and }\\
\text{locally finite } \quad &\text{ if } \#(Y \cap \La) < \infty
\text{ for all } \La \in \B^2_b.
\end{array}
\] 
The configuration space  $\Y$ of particles is defined as the set of all locally finite subsets of $\RS$. The elements of $\Y$ are called configurations of particles.  For $Y \in \Y$ and $A \in \B^2 \otimes \F_S$ let  
\[
\begin{array}{ll}
Y_{A} := Y \cap A 
&\text{(restriction of $Y$ to $A$)}, \\
\Y_{A} := \{ Y \in \Y: Y \subset A \} 
&\text{(set of all configurations in $A$),  and }\\
N_{A}(Y) := \# Y_{A}  
&\text{(number of particles of $Y$ in $A$).} 
\end{array}
\]
The counting variables $(N_{A})_{A \in \B^2 \otimes \F_S}$ generate a $\si$-algebra 
on $\Y$, which will be denoted by  $\F_{\Y}$. For $\La \in \B^2$ let 
$\F'_{\Y,\La}$ be the $\si$-algebra on $\Y_{\La}$ obtained by restricting
$\F_{\Y}$  to $\Y_{\La}$, and let $\F_{\Y,\La} := e_{\La}^{-1} \F'_{\Y,\La}$ be the 
$\si$-algebra on $\Y$ obtained from $\F'_{\Y,\La}$ by the restriction mapping 
$e_{\La}: \Y \to \Y_{\La}, Y \mapsto Y_{\La}$. The tail $\si$-algebra or
$\si$-algebra of the events far from the origin is defined by 
$\F_{\Y,\infty} \, := \, \bigcap_{n \ge 1} \F_{\Y,\Lan^c}$.
For configurations $Y,\bY \in \Y$ let  $Y\bY := Y \cup \bY$.
Let $\nu$ be the distribution of the Poisson point process on 
$(\Y,\F_{\Y})$, i.e. 
\[
\int \nu(dY) f(Y)  \, = \, e^{- \leb(\La)} \, \sum_{ k \ge 0} \,  
 \frac{1}{k!} \, \int_{{\La}^k} dy_1 \ldots dy_k \, 
 \,  f(\{y_i : 1 \le i  \le k\} ) 
\]
for any $\F_{\Y,\La}$-measurable function 
$f: \Y \to \R_+$, where $\La \in \Bo_b^2$.  
For $\La \in  \Bo^2_b$ and $\bY \in \Y$ 
let $\nu_{\La}(.|\bY)$ be the distribution of the Poisson point 
process in $\La$ with boundary condition $\bY$, i.e.
\begin{equation*} 
\int \nu_{\La}(dY|\bY) f(Y) \;
= \; \int \nu (dY) f(Y_{\La}\bY_{\La^c})
\end{equation*}
for any $\F_{\Y}$-measurable function $f: \Y \to \R_+$.
It is easy to see that $\nu_{\La}$ is a stochastic kernel 
from $(\Y,\F_{\Y,{\La}^c})$ to  $(\Y,\F_{\Y})$.\\

The configuration space of bonds is the set of  all locally 
finite bond sets:
\[
\E \, := \,  \{B \subset E(\RS): \, 
\#\{ yy' \in B : yy' \subset \La \times S \} < \infty \, 
\text{ for all } \La \in \B_b^2  \}.
\]
On $\E$ the $\si$-algebra $\F_{\E}$ is defined to be generated by the 
counting variables 
$N_{E}: \E \to \N,  B \mapsto \#(E \cap B)$ $\;(E \in \F_{E(\RS)})$.  
For a countable set $E \in \E$ one can also consider the 
Bernoulli-$\sigma$-algebra $\Bo_E$ on $\E_{E} := \p(E) \subset \E$,  
which is defined to be generated by the family of sets 
$(\{B \subset E: e \in B\})_{e \in E}$. Given a family 
$(p_e)_{e \in E}$ of reals in $[0,1]$ the Bernoulli
measure on $(\E_E,\B_E)$ is defined as the unique probability measure  
for that the events $( \{B \subset E : e \in B \} )_{e \in E}$ 
are independent with probabilities $(p_e)_{e \in E}$.
It is easy to check that the inclusion $(\E_E,\B_E) \to (\E, \F_{\E})$ is 
measurable. Thus any probability measure on  
$(\E_E,\B_E)$ can trivially be extended to $(\E, \F_{\E})$.


\subsection{Gibbs measures} \label{secgibbs}

Let $U:(\RS)^2 \to \baR$ be a potential and $z >0$ an 
activity parameter. For finite configurations $Y,Y' \in \Y$ 
we consider the energy terms 
\[
H^U (Y) \, := \sum_{y_1y_2 \in E(Y)} U(y_1,y_2) \quad \text{ and } \quad  
W^U (Y,Y') \, := \sum_{y_1 \in Y} \sum_{y_2 \in Y'} U(y_1,y_2).
\]
The last definition can be extended to locally finite configurations  $Y'$ 
whenever  $W^U (Y,Y'_{\La})$ converges as  $\La \uparrow \R^2$ 
through the net  $\B^2_b$. The Hamiltonian of a configuration $Y \in \Y$ in 
$\La \in \B^2_b$ is given by 
\[ 
\begin{split}
H^U_{\La}&(Y) \, := \, 
H^U (Y_{\La}) +  W^U(Y_{\La},Y_{\La^c}) \, =  
\sum_{y_1y_2 \in E_{\La} (Y)}  U(y_1,y_2),\\ 
& \text{ where } \quad  
E_\La (Y) := \{ y_1y_2 \in E(Y):  y_1 y_2 \cap \La \ne \emptyset\}.
\end{split}
\]
The integral 
\[ 
Z^{U,z}_{\La}(\bY) \, 
:= \, \int \nu_{\La}(dY|\bY) \, e^{-H^U_{\La}(Y)}z^{\#Y_\La}
\]
is called the partition function in  $\La \in \B^2_b$ for the boundary condition 
$\bY_{\La^c} \in \Y$. In order to ensure that the above objects are 
well defined and the partition function is finite and positive 
we need the following definition: 
\begin{defi} \label{defzul}
A triple $(U,z,\YY)$ consisting of a potential $U:(\RS)^2 \to \baR$, 
an activity parameter $z >0$, and a set of boundary conditions  
$\YY \in \F_{\Y,\infty}$ is called admissible if for all $\bY \in \YY$ and 
$\La \in \B^2_b$ the following holds: 
$W^U(\bY_{\La},\bY_{\La^c})$ has a well defined value in $\baR$, and 
the partition function  $Z^{U,z}_{\La}(\bY)$ is finite.
\end{defi}
If $(U,z,\YY)$ is admissible, $\La \in \B^2_b$, and $\bY \in \YY$, 
then $W^U(Y_{\La},\bY_{\La^c}) \in \baR$ is well defined for every $Y \in \Y$, 
because $Y_{\La}\bY_{\La^c} \in \YY$. As a consequence the partition function 
$Z^{U,z}_{\La}(\bY)$ is well defined. Furthermore by definition it is finite, and 
by considering the empty configuration one can show that it is positive. 
The conditional Gibbs distribution $\gmu^{U,z}_{\La}(.|\bY)$ in $\La \in \B^2_b$ 
with boundary condition $\bY \in \YY$ is thus well defined by 
\[
\gmu^{U,z}_{\La}(A|\bY) \, := \, 
\frac 1 {Z^{U,z}_{\La}(\bY)}  \int \nu_{\La}(dY|\bY) \, 
e^{-H^U_{\La}(Y)} z^{ \#Y_\La} 1_A(Y) \quad 
\text{ for } \quad  A \in \F_{\Y}.
\]
$\gmu^{U,z}_{\La}$ is a probability kernel from $(\YY,\F_{\YY,\La^c})$ to
$(\Y,\F_{\Y})$. Let 
\[ 
\begin{split}
\G_{\YY}(U,z) \, := \, \{ \mu \in &\prob(\Y,\F_{\Y}) : \, \mu(\YY) = 1 
\quad \text{ and } \\
&\mu(A|\F_{\Y,\La^c}) = 
 \gmu^{U,z}_{\La}(A|.) \text{ $\mu$-a.s. } \alle A \in \F_{\Y},
\La \in \B^2_b \}
\end{split}
\]
be the set of all Gibbs measures corresponding to $U$ and $z$ with whole 
weight on boundary conditions in $\Y_0$. It is easy to see that for any 
probability measure $\mu \in \prob(\Y,\F_{\Y})$ such that $\mu(\YY)=1$ 
we have the equivalence 
\[ 
\mu \in \G_{\YY}(U,z) \quad \gdw \quad (\mu \otimes \gmu^{U,z}_{\La} = \mu
 \alle \La \in \B^2_b).
\] 
So for every  $\mu \in \G_{\YY}(U,z)$, $f: \Y \to \R_+$ measurable 
and $\La \in \B^2_b$ we have 
\begin{equation} \label{gibbsaequ}
\int \mu(dY) \, f(Y)\, = \, 
\int \mu(d\bY) \int \gmu_{\La}^{U,z} (dY|\bY) \, f(Y). 
\end{equation}
We note that the hard core $\KU$ of a potential $U$ models the property 
that particles are not allowed to get too close to each other, i.e.
for admissible $(U,z,\YY)$ and $\mu \in \G_{\YY}(U,z)$ we have
\begin{equation} \label{hardcore}
\mu(\{Y \in \Y: \exists y,y' \in Y: y\neq y',(y,y') \in \KU \})
\, = \, 0. 
\end{equation}
This is because for every $n \in \N$ and every boundary condition 
$\bY \in \YY$ we have  
\[
\gmu_{\Lan}^{U,z} (\{Y \in \Y: \exists y,y' \in Y_{\Lan}: y\neq y',(y,y') 
\in \KU \} | \bY) \, = \, 0, 
\]
as on the event considered in the last line the Hamiltonian 
$H_{\Lan}^U(Y|\bY)$ is infinite. Therefore \eqref{hardcore} follows 
by using \eqref{gibbsaequ} and taking $n \to \infty$.\\

For admissible $(U,z,\YY)$ and a Gibbs measure $\mu \in \G_{\YY}(U,z)$ 
we define the correlation function  $\rho^{U,\mu}$ by 
\[
\rho^{U,\mu}(Y) \, = \, e^{-H^U(Y)} \int \mu (d\bar{Y}) \,  
     e^{- W^U (Y,\bar{Y})}
\]
for any finite configuration $Y \in \Y$. If there is a $\xi = \xi(U,z,\YY) \ge 0$ 
such that
\begin{equation} \label{ruellebd}
\rho^{U,\mu}(Y) \, \le \, \xi^{\#Y} \quad \text{ for all finite }
Y \in \Y \text{ and all } \mu \in \G_{\YY}(U,z),  
\end{equation}
then we call $\xi$ a Ruelle bound for  $(U,z,\YY)$. Actually we need 
this bound on the correlation function in the following way:
\begin{lem} \label{lekorab}
Let $(U,z,\YY)$ be admissible with Ruelle bound $\xi$.
For every Gibbs measure $\mu \in \G_{\YY}(U,z)$ and every measurable function 
$f: (\RS)^m \to \R_+$, $m \in \N$, we have 
\begin{equation} \label{korab}
\int \mu(dY) \sideset{}{^{\neq}}\sum_{y_1,\ldots,y_m \in Y} 
f(y_1,\ldots,y_m) \,
\le \, (z\xi)^m  \int dy_1 \ldots dy_m  \, f(y_1,\ldots,y_m).
\end{equation}
\end{lem}
We use $\Sigma^{\neq}$ as a shorthand notation for a multiple sum such that
the summation indices are assumed to be pairwise distinct.


\subsection{Superstability and admissibility} \label{secadm} 

Now we will discuss some conditions on potentials that imply that $(U,z,\YY)$ 
is admissible and has a Ruelle bound whenever the set of boundary conditions
$\YY$ is suitably chosen. Apart from purely repulsive (i.e nonnegative)
potentials such as the one considered in Theorem~\ref{pottstype}
we also want to consider superstable potentials in the sense of 
Ruelle, see \cite{Ru1}. Therefore let  
$\Gar  :=   r + [-\frac 1 2 , \frac 1 2 [^2 \subset \R^2$  
be the unit square centred at  $r \in \Z^2$
and let  $\Z^2(Y)  :=  \{ r \in \Z^2: N_{\Gar}(Y) > 0 \}$
be the minimal set of lattice points such that the corresponding squares cover the configuration $Y$. A potential $U:(\RS)^2 \to \baR$ 
is called superstable if there are real constants  $a > 0$ and $b \ge 0$ 
such that
\[ 
H^U (Y) \, \ge \,\sum_{r \in \Z^2(Y)} \,[a N_{\Gar}(Y)^2 - b N_{\Gar}(Y)]
\]
for all finite configurations  $Y \in \Y$.
$U$ is called lower regular if there is a decreasing function $\Psi: \N \to \R_+$ 
with $\sum \limits_{r \in \Z^2} \Psi(|r|) <  \infty$ such that
\[
W^U (Y,Y') \, \ge \, -\sum_{r \in \Z^2(Y)} 
\sum_{s \in \Z^2(Y')} \, \Psi(|r-s|) \, 
[\frac{1}{2} N_{\Gar}(Y)^2 + \frac{1}{2} N_{\Gas}(Y')^2] 
\]
for all finite configurations $Y,Y' \in \Y$. So superstability and 
lower regularity give lower bounds on energies 
in terms of particle densities. In order to control these densities
a configuration  $Y \in \Y$ is said to be tempered if 
\[
\bs(Y) \, := \, \sup_{n \in \N} s_n(Y) \, < \, \infty, \; \text{ where } \; 
s_n(Y) \, := \, \frac{1}{(2n+1)^2} \sum_{r \in \Z^2 \cap \La_{n+1/2}} N_{\Gar}^2(Y).
\]
By $\Y_t$ we denote the set of all tempered configurations. We note that 
$\Y_t \in \F_{\Y,\infty}$.  
\begin{lem} \label{lemsup}
Let $z > 0$ and $U:(\RS)^2 \to \baR$ be a potential function.
\begin{enumerate}
\item[(a)]
If $U \ge 0$, then  $(U,z,\Y)$ is admissible 
with Ruelle bound  $\xi := 1$.
\item[(b)]
If $U$ is superstable and lower regular then $(U,z,\Y_t)$ is admissible and
admits a Ruelle bound. 
\end{enumerate}
\end{lem}
The first assertion is a straightforward consequence of the fact that
all energy terms are nonnegative. For the second assertion see \cite{Ru1}.


\subsection{Conservation of translational symmetry} \label{secsym}

We want to establish the conservation of $\tauh$-translational 
symmetry for a given admissible triple  $(U,z,\YY)$ and a translation 
$\tauh \in \R^2$. It suffices to show  that for every $\delta >0$ and every cylinder event 
$D \in \F_{\Y,\La_{m}}$  $(m \in \N)$ there is a natural  $n \ge m$ such that 
we have
\begin{equation} \label{georgiianalog}
\gmu^{U,z}_{\Lan}(D + \tauh|\bY) \, + \, \gmu^{U,z}_{\Lan}(D- \tauh|\bY) \, 
\ge \, 2 \gmu^{U,z}_{\Lan}(D|\bY) \, - \, \delta \quad \text{ for all } \bY \in \YY.
\end{equation} 
Indeed, let  $\mu \in \G_{\YY}(U,z)$, then integrating  \eqref{georgiianalog}
with respect to $\mu$ and applying \eqref{gibbsaequ} gives  
$\mu(D + \tauh)  +   \mu(D - \tauh) \ge  2 \mu(D)  -  \delta$
for all $\delta > 0$ and $D \in \F_{\Y,\La_{m}}$  $(m \in \N)$. Letting $\delta \to 0$, 
Lemma~\ref{lekrit} shows the invariance of $\mu$ under the translation  by $\tauh$, 
because the cylinder events form an algebra which generates the $\si$-algebra 
$\F_{\Y}$. For the proof of Theorem~\ref{sym} our goal 
thus will be to establish an inequality similar to \eqref{georgiianalog}.
We further note that the group $\R \ah$ is generated by the set 
$\{\tau \ah: \tau \in [0,1/2]\}$. Thus we only have to consider translations of this 
special form in order to establish the $\ah$-invariance of a set of Gibbs measures.


\subsection{Concerning measurability} \label{secmeas}

We will consider various types of random objects, all of which have
to be shown to be measurable with respect to the considered $\si$-algebras. 
However, we will not prove measurability of every such object in detail.
Instead,  we will now give a list of operations that preserve
measurability.
\begin{lem} \label{lemeas}
Let $Y,Y' \in \Y$, $B,B' \in \E$, $y \in \RS$, $x \in \R^2$, and  
$p \in \Omega$ be variables, 
where $(\Omega,\F)$ is a measurable space. 
Let $f: \Omega \times (\RS) \to \R$ and 
$g: \Omega \times E(\RS) \to \R$ be measurable. 
Then the following functions of the given variables are measurable 
with respect to the considered $\si$-algebras: 
\begin{align}
&\sum_{y' \in Y} f(p,y'), \quad Y \cap Y', \quad Y \cup Y', 
 \quad Y \weg Y', \quad Y+x,\label{measY}\\
&\sum_{b' \in B} g(p,b'), \quad B \cap B', \quad B \cup B', 
 \quad B \weg B',  \quad  B+x, \label{measB}\\
& \inf_{y' \in Y} f(p,y'), \quad 
C_{Y,B}(y),\quad E(Y) \label{measinf}.
\end{align}
\end{lem}
Using this lemma and well known theorems, such as the measurability part of 
Fubini's theorem, we can check the measurability of 
all objects considered. 

%
\newpage

\section{Proof of the lemmas from Sections~\ref{secresult} and 
\ref{secsetting}} \label{secleset}

\subsection{Conservation of symmetries: Lemma~\ref{lekrit}} 

We first note that  $\mu \circ \tau^{-1} = \mu$ easily implies 
\eqref{kritA}, where we indeed have equality. 
For the other implication let us assume \eqref{kritA}. 
By the monotone class theorem this inequality immediately extends
to all $A \in \F$. Thus for all $D \in \F$ and $k\in\Z$ 
\[
\mu(\tau^{k+1} D) \, + \, \mu(\tau^{k-1} D) \, \ge \, 2 \mu(\tau^k D),
\]
i.e. the sequence $(\mu(\tau^k D))_{k \in \Z}$ is convex. But 
$\mu$ is a probability measure, so the sequence is bounded, and 
thus it has to be constant. In particular we get 
$\mu(\tau^{-1}D) = \mu(D)$. As $D \in \F$ 
was arbitrary the result follows.


\subsection{Measurable enlargements: 
Lemmas~\ref{lebedeperw} and \ref{leballs}}

Let $A$ be as described in Lemma~\ref{lebedeperw} and let $\ep >0$. 
By Fubini's theorem  the function $f: (\RS)^2 \to \R_+$, 
$f(y_1,y_2) \, := \la(A(y_1,y_2,\ah) \, \cap\;  ]-\ep,\ep[\, )$
is measurable and by \eqref{bedeperw} 
for all $y_1,y_2\in \RS$ we have 
\[
f(y_1,y_2) > 0 \quad 
\gdw  \quad A(y_1,y_2,\ah) \, \cap\;  ]-\ep,\ep[ \;\neq \emptyset \quad 
\gdw \quad (y_1,y_2) \in A_{\ep,\ah}.
\]
This shows $A_{\ep,\ah} = \{ f > 0 \}$ to be measurable. 
The second statement of Lemma~\ref{lebedeperw} is an immediate 
consequence  of the fact that a Borel set containing a nonempty 
open set has positive Lebesgue-measure. \\

For the proof of Lemma~\ref{leballs} (a) let all $\ah$-cross 
sections of $\Ko$ be open intervals.  Then
the $\ah$-cross sections $(\Ko_{\ep,\ah} \weg \Ko)(y_1,y_2,\ah)$ are either empty or the union of two intervals of length $\ep$. Furthermore, as $\Ko$ is of bounded range there is a real $r>0$ such that for every 
$(y_1,y_2) \in \Ko_{\ep,\ah}$ we have $|y_1-y_2|\le r$.
The supremum in (a) can thus be estimated by  
\[
\sup_{r_1,r_1' \in \R} \sup_{\si_1,\si_2 \in S} 
\int dr'_2 1_{\{|r'_1-r'_2| \le r\}} \int dr_2  
1_{\{r_2 \in (\Ko_{\ep,\ah} \weg \Ko)(r_1,r'_1,\si_1,r'_2,\si_2)\}} \, \le \, 2r \cdot 2 \ep.
\] 
For part (b) let $\Ko$ be a standard set consisting of discs and 
let $r >0$ be a bound for the corresponding family $(r_{\si_1 \si_2})_{\si_1,\si_2 \in S}$. The supremum in (b) can be estimated by  
\[
\sup_{x_1 \in \R^2} \sup_{\si_1,\si_2 \in S} 
\int 1_{\{r_{\si_1\si_2} < |x_1-x_2|_h  \le r_{\si_1\si_2} + \ep\}} dx_2 \; 
\le \; \leb(\{ r < |.|_h< r+\ep \}). 
\] 


\subsection{Smooth or continuous potentials: 
Lemma~\ref{lepotglatt}}

We set $\Ko := \KU_{+\de}$, $\Kop := \Ko_{+\ep}$, $\Kopp =\Kop_{+\ep}$, 
where $\ep,\de > 0$ are so small that 
\[
c \, :=  \, 1/(z\xi) - \sup_{y_1 \in \RS} 
             \int 1_{\Kopp \weg \KU} (y_1,y_2) \, dy_2 \, > \, 0,
\]
where $\xi$ is a Ruelle bound for $(U,z,\YY)$. (This is possible by Lemma~\ref{leballs}.) In case (a) we are done
setting $\psi := \psi _\de$, $\bU := U$ and $u :=0$. 
In case (b1) let $U_1 := U$ and in case (b2) let $U_1 :=  \tU$.
Without loss of generality we may assume that  
$R \ge 1$ and $\Ko \subset \tilde{K}$, and furthermore
\[
\sup_{y_1 \in \RS} \int  2 \tU(y_1,y_2) |y_1-y_2|^2 
1_{\tilde{K}^c}(y_1,y_2) \, dy_2 \, < \, \frac{c}{2} 
\]
in case (b2). In both cases  $U_1$ serves as an $\ah$-smooth approximation of $U$ on 
$\tilde{K}^c$. We note that $U_1$ is bounded and has $\psi_{\delta}$-dominated 
$\ah$-derivatives on $\tilde{K}^c$, which also implies that 
$\partial_{\ah}^2 U_1$ and $\partial_{\ah} U_1$ are bounded on  $\tilde{K}^c$. Let 
\[
C := \{(y_1,y_2) \in (\RS)^2: |y_1 -y_2| \le R + 1\}  \weg \Ko.
\]
For  $\delta'> 0$ let $f_{\delta'}: \R \to \R_+$ be a symmetric smooth probability 
density with support in   $]-\delta',\delta'[$, e.g. 
$f_{\delta'}(t) := \frac 1{c_{\delta'}} 1_{]-\delta',\delta'[}(t) e^{-(1-t^2/\delta'^2)^{-1}}$, where $c_{\delta'}$ is a normalizing 
constant. Then 
\[
U_2(x_1,\si_1,x_2,\si_2) \, 
:= \, \int dt \,  f_{\delta'}(t)  U(x_1,\si_1,x_2 - t\ah,\si_2)
\]
is an $\ah$-smooth approximation of $U$ on $C$. If  $\delta'$ is small enough, 
then
\[
|U_2(y_1,y_2) - U(y_1,y_2)| \, <  \, \frac{c}{16(R+1)^2} 
\quad \text{ for  } (y_1,y_2) \in C
\]
by the $\ah$-equicontinuity of $U$. 
Let  $g: (\RS)^2 \to [0,1]$ be an $\ah$-smooth function with 
$g(y_1,y_2) = 0$ for $|y_1-y_2| \le R$, $g(y_1,y_2) = 1$  for $|y_1-y_2| \ge R+1$ 
and such that the  $\ah$-derivatives $\partial_{\ah} g$ and $\partial^2_{\ah} g$
are bounded. Now we can define $\bU, u: \Ko^c \to \R$ by  
$\bU  :=  (1 - g) (U_2 + c') + g U_1$  and $u  :=  \bU - U$. 
It is easy to verify that the constructed objects have all the properties 
described in Definition~\ref{defGapprox} in both cases (b1) and (b2).


\subsection{Potts type potentials: Lemma~\ref{lepotpotts}}

We first consider a well behaved function $\phi: \R_{+} \to \baR$ 
(with respect to given reals $0 \le r_0 < \ldots < r_n$, $n \ge 0$)
and show how to decompose $\phi$ into a continuous part $\bar{\phi}$ 
and a small part $\varphi$. 
For $s,\ep > 0, m\in \R$ we define $h_{s,m,\ep}: \R_+ \to \R$ such that 
the graph of $h_{s,m,\ep}$ looks like $\bigwedge$, where $(s,m)$ is the 
topmost point and the angle is determined by $\ep$, i.e. 
$h_{s,m,\ep}(r) := m- (m/ \ep)|r-s|$. Defining
\[
\bar{\phi} := \phi \vee \bigvee_{i=1}^n h_{r_i,m_i,\ep}, \quad \text{ where } \; 
m_i := \Big( \phi(r_i) \vee \lim_{r \to r_i+} \phi(r) \vee \lim_{r \to r_i-}\phi(r)\Big) + 1, 
\]
for a given $\ep>0$, we see that $\bar{\phi}$ is continuous on $]r_0,\infty[$. 
Furthermore $\varphi := \bar{\phi} - \phi$ satisfies $\varphi \ge 0$, 
$\int (\varphi(|y|_h) \wedge 1)|y|^2 1_{\{|y|_h > r_0 \}} dy < \infty$, and 
$\int (\varphi(|y|_h) \wedge 1)  1_{\{|y|_h > r_0 \}} dy$ is 
arbitrarily small if only $\ep >0$ is chosen small enough.\\

Now let $U$ be a Potts type potential corresponding to a norm $|.|_h$ 
and a family of well behaved interactions $(\phi_{\si_1 \si_2})_{\si_1,\si_2 \in S}$, 
where $S$ is a finite spin space. As above we decompose every 
$\phi_{\si_1 \si_2}$ into a continuous part $\bar{\phi}_{\si_1 \si_2}$ and 
a small part $\varphi_{\si_1 \si_2}$, where the $\ep>0$ entering 
the above construction is chosen sufficiently small. Now let 
$U_c(y_1,y_2) := \bar{\phi}_{\si_1 \si_2}(|x_1-x_2|_h)$ and $u_c := U_c - U$. 
We observe that $U_c$ is of the form described in Lemma~\ref{lepotglatt} (b1): 
We simply choose $\psi = 0$ and $\tilde{K}$ so big such that $U_c= 0$ on 
$\tilde{K}^c$. We note that $U_c$ is bounded and $\ah$-equicontinuous in 
$\KU_{+\de}$ for every $\de >0$. As in the proof of 
Lemma~\ref{lepotglatt} we thus find a decomposition of $U_c$ into 
suitable potentials $\bU$ and $u$. Then $\bU$ and $u + u_c$ give a 
decomposition of $U$ into a smooth part and a small part as required.


\subsection{Property of the Ruelle bound: Lemma~\ref{lekorab}} 

For every $n \in \N$, every measurable $g: \Y_{\La_n} \to \R_+$ 
and every $\bY \in \YY$ we have 
\begin{displaymath}
\begin{split}
\int &\nu_{\Lan}(dY|\bY)  \sideset{}{^{\neq}}\sum_{y_1,\ldots,y_m \in 
    Y_{\La_n}} \! \! f(y_1,\ldots,y_m) \, g(Y) \\ 
&= \, \int_{{\Lan}^m} dy_1 \ldots dy_m 
  \, f(y_1,\ldots,y_m) \,
  \int \nu_{\Lan}(dY'|\bY)\, g( \{y_1, \ldots, y_m \} Y').
\end{split}
\end{displaymath}
Combining this with \eqref{gibbsaequ}, the definition of the conditional 
Gibbs distribution and the definition of the correlation function we get 
\begin{displaymath}
\begin{split}
&\int \mu(dY) \sideset{}{^{\neq}}\sum_{y_1,\ldots, y_m \in Y_{\Lan}} \!\!
 f(y_1,\ldots ,y_m)\\
&= \,  \int \mu(d\bY) \, \frac{1}{Z^{U,z}_{\La_{n}}(\bY)}
   \, \int \nu_{\Lan} (dY|\bY) \sideset{}{^{\neq}} 
   \sum_{y_1,\ldots ,y_m \in Y_{\Lan}}  \!\! f(y_1,\ldots ,y_m) \, 
    e^{- H^U_{\La_{n}}(Y)}z^{ \# Y_{\Lan}}\\   
&= \, \int_{{\Lan}^m} dy_1 \ldots  dy_m \, f(y_1,\ldots ,y_m)\, z^m 
   \, \rho^{U,\mu} (\{y_1,\ldots ,y_m\}).  
\end{split} 
\end{displaymath} 
Now we use  \eqref{ruellebd} to estimate the correlation function by 
the Ruelle bound  $\xi$.  Letting  $n \to \infty$ the assertion follows
from the monotone limit theorem.


\subsection{Measurability: Lemma~\ref{lemeas}}

Details concerning measurability of functions of point processes 
can be found in \cite{DV} or \cite{MKM}, for example. 
The first part of \eqref{measY} is the 
measurability part of Campbell's theorem.
For the rest of \eqref{measY} it suffices to observe that  we have  
$N_{A}(Y \cap Y') = \sum \limits_{y \in Y} 
\sum \limits_{y' \in Y'} 1_{\{y = y' \in A\}}$, 
$N_A(Y \weg Y') =  N_A(Y) -  N_A(Y \cap Y')$, 
$N_A(Y \cup Y') =  N_A(Y) +  N_A(Y' \weg Y)$ and 
$N_A(Y+x) = \sum \limits_{y \in Y} 1_A(y+x)$ for all
$A \in \Bo^2_b \otimes \F_S$. 
\eqref{measB} can be proved similarly. For \eqref{measinf}
we note that
$\inf \limits_{y' \in Y} f_1(p,y') < c$  
$\gdw \sum \limits_{y' \in Y} 1_{\{f_1(p,y')<c\}} \ge 1$ 
for all  $c \in \R$, 
$N_A(C_{Y,B}(y))  =  \sum \limits_{y' \in Y} 
 1_{\{y' \in C_{Y,B}(x), y' \in A\}}$ for all $A \in \B^2_b \otimes \F_S$,
$y' \in C_{Y,B}(y)$ $\gdw \sum \limits_{m \ge 0} 
\sum \limits_{y_0,\ldots ,y_m \in Y} 1_{\{y=y_0,y'=y_m\}}
\prod \limits_{i=1}^m 1_{\{ y_i y_{i+1} \in B  \}} \ge 1$ for 
all $y' \in \RS$ and 
$N_L(E(Y))  =  \frac 1 2 \sum_{y_1,y_2 \in Y} 1_{\{y_1y_2\in L\}}$
for all $L \in \F_{E(\RS)}$.
Using these relations, the measurability of the terms in \eqref{measinf} follows easily. Note that we made repeated use of the fact that the 
diagonal is measurable in $S \times S$.  

%

\newpage

\section{Proof of Theorem~\ref{sym}: Main steps} \label{secproofsym}

\subsection{Basic constants} \label{constants}

Let $(U,z,\YY)$ be admissible with Ruelle bound $\xi$, where 
$U:(\RS)^2 \to \baR$ is an $\ah$-invariant, $\ah$-smoothly approximable standard potential. There is an $\ah$-invariant standard set  $\Ko \supset \KU$ with $\ah$-enlargements $\Kop$, $\Kopp$,
a bpsi-function $\psi$, and measurable symmetric
$\ah$-invariant functions  $\bU:  (\RS)^2   \to \baR$ and 
$u: (\RS)^2 \to \R$ such that $U =  \bU - u$ and 
$u  \ge  0$  on $(\RS)^2$, $u = 0$  on $\Ko$, 
$\bU$ has $\psi$-dominated $\ah$-derivatives on $\Ko^c$,  and 
$\tiu  =  1 - e^{-u}$ satisfies 
\begin{equation} \label{dec}
\begin{split}
&\cu \, := \, \sup_{y_1 \in \RS} \int \tiu (y_1,y_2) |y_1-y_2|^2 dy_2 
\,  < \, \infty \quad \text{ and } \\
&\ex \, := \,  \sup_{y_1 \in \RS}  \int  (1_{\Kopp \weg \KU} + \tiu) 
     (y_1,y_2) dy_2 \, < \, \frac{1}{z\xi}. 
\end{split}
\end{equation}
Note that above we defined $\bU$ and $u$ 
also on $\Ko$. 
By symmetry  we may suppose that the direction of the translations is  
$\ah = \einh := (1,0)$, and w.l.o.g.  we may assume that 
\[
K \supset \{(x_1,\si_1,x_2,\si_2) \in (\RS)^2: x_1 = x_2\}.
\]
Let $\fc:(\RS)^2 \to [0,1]$ be 
a measurable function  such that 
\[
\fc = 0 \, \text{ on } \, \Ko, \; \fc = 1 \, \text{ on } \, (\Kopp)^c,
 \; \fc \text{ is $\einh$-smooth, and $\partial_{\einh} \fc$ is bounded},
\]  
where $\partial_{\einh} \fc$ is the $\einh$-derivative with respect to the
second spatial component. For the construction of such a function 
we introduce $\ft: (\RS)^2 \to \R$, $\ft:= 1_{(\Kop)^c}$,
and choose an infinitely often differentiable function 
$f_{\ep}: \R \to \R_+$ 
which is a probability density with support in  $]-\ep,\ep[$. Then 
the function
$\fc(y',y)  :=  \int dt    \ft(y',y- t \einh)  f_{\ep}(t)$
has the desired properties if $\ep>0$ is chosen small enough. Furthermore we need the following constants: 
\begin{equation} \begin{split} \label{bpsi}  
\cpsi \, &:= \, \|\psi\|\, \vee \sup_{y_1 \in \RS} 
\int dy_2  \, \psi(y_1,y_2)(|y_1-y_2|^2 \vee 1) , \\  
\cK \, & := \, \sup\{|y_1-y_2|: (y_1,y_2)  \in \Kopp \}, \quad \text{ and }
\quad \cf \, := \, \| \partial_{\einh} \fc \|. 
\end{split} \end{equation}
These constants are finite as $\psi$ is a  bpsi-function, $\Kopp$ has 
bounded range, and $\partial_{\einh} \fc$ is bounded. On $\RS$ 
we consider the partial order $\lear$ defined by 
\[
 (r_1,r_2,\si) \, \lear \, (r_1',r_2',\si') \quad 
:\Leftrightarrow \quad r_1 \le r_1', r_2 = r_2',\si=\si'.
\]
In order to show the conservation of $\einh$-translational symmetry 
we fix a Gibbs measure $\mu \in \G_{\YY}(U,z)$, a cylinder event 
$D \in \F_{\Y,\La_{n'-1}}$, where $n'\in \N$, a real $\de \in]0,1/2[$,
and a translation distance parameter  $\tau  \in [0,1/2]$, see
subsection \ref{secsym}. We will ignore dependence on any of the 
above parameteres in our notations.


\subsection{Decomposition of $\mu$ and the bond process}

We consider the bond set $\en(Y)  :=   E_{\Lan}(Y)  =  
\{ y_1y_2 \in E(Y):  y_1 y_2 \cap \Lan \ne \emptyset \}$
for $n \in \N$ and $Y \in \Y$. 
On $(\E_{\en(Y)},\B_{\en(Y)})$ we introduce the Bernoulli measure 
$\pin(.|Y)$ with bond probabilities 
\[
(\tiu(b))_{b \in \en(Y)} \quad \text{ where } \quad  \tiu(b) \, := \,  1-e^{-u(b)}, 
\]
using the shorthand notation $u(y_1y_2) := u(y_1,y_2)$ for 
$y_1,y_2 \in \RS$. 
We  note that   $0 \le \tiu(b) < 1$ for all $b \in \en(Y)$ as 
$0 \le u < \infty$. As remarked earlier $\pin(.|Y)$ can be extended to 
a probability measure on $(\E, \F_{\E})$. For all $D \in \F_{\E}$ $~\pin(D|.)$ 
is $\F_{\Y}$-measurable, so $\pin$ is a probability kernel from $(\Y,\F_\Y)$ 
to $(\E,\F_\E)$. 
\begin{lem} \label{umord}
Let $n \in \N$. We have  
\[   
\mu \otimes \nu_{\Lan}(\Gnb) = 1 \;  \text{ and } \; \mu(\Gnb) = 1 \; 
\text{ for } \; \Gnb  :=  \Big\{Y \in \Y: \sum \limits_{b \in \en(Y)} 
 \tiu(b)   < \infty \Big\}.
\]
\end{lem} 
For $Y \in \Gnb$ every bond set  
is finite $\pin(.|Y)$-a.s. by Borel-Cantelli, so 
$\pin( . |Y) \ll \pinh(.|Y)$, where $\pinh(.|Y)$ denotes 
the counting measure on  $(\E_{\en(Y)},\B_{\en(Y)})$ 
concentrated on finite bond sets.  Again, $\pinh$ can be considered 
as a probability kernel from $(\Y,\F_\Y)$ to $(\E,\F_\E)$. 
We have  
\[
\frac{d\pin(.|Y)}{d\pinh(.|Y)} (B) \, 
= \, \prod_{b \in B} \tiu(b)  \prod_{b \in \en(Y) \weg B} (1 - \tiu(b)) \,
= \, e^{-H^{{u}}_{\Lan}(Y)}  \prod_{b \in B} (e^{u(b)}-1),  
\]
so for every $Y \in \Gnb$ the Hamiltonian $H^u_{\Lan}(Y)$ is finite, 
and thus the decomposition of the potential gives a corresponding
decomposition of the Hamiltonian 
\[
H^U_{\Lan}(Y) \, = \,  H^{\bar{U}}_{\Lan}(Y) - 
H^{{u}}_{\Lan}(Y).
\]
Using \eqref{gibbsaequ} we conclude that for every 
$\F_\Y \otimes \F_{\E}$-measurable function $f \ge 0$  
\begin{equation}  \label{desint} 
\begin{split} 
\int d\mu \otimes \pin \, f \;  
&= \;  \int \mu(d\bY) \frac 1  {Z_{\Lan}^{U,z}(\bar{Y})} \int \nu_{\Lan} \otimes \pinh (dY,dB |\bY)  \\
&\hspace{ 2 cm} z^{\#Y_{\Lan}} e^{-H^{\bar{U}}_{\Lan}(Y)} 
 \prod_{b \in B} (e^{u(b)}-1)  f(Y,B). 
\end{split} \end{equation}
Here by  Lemma~\ref{umord} on both sides we have  $Y \in \Gnb$ with 
probability one, thus the equality follows from the above decomposition. 
If $f$ does not depend on $B$ at all, 
the integral on the left hand side of \eqref{desint} is just the $\mu$-expectation 
of $f$, as $\pin$  is a probability kernel, and from the right hand side
we learn that the perturbation $u$ of the $\einh$-smooth potential $\bU$ 
can be encoded in a bond process $B$ such that the perturbation affects
only those pairs of particles with $y_1y_2 \in B$.  


\subsection{Generalized translation and good configurations} \label{tgt} 

For integers $n,R$ such that $n > R \ge n'$ we define the functions
$q: \R_+ \to \R$, $Q: \R_+ \to \R$, $r: \R \times \R_+ \to \R$  and  
$\tn: \R \to \R$ by
\begin{displaymath}
\begin{split}
q(s) \, &:= \, \frac{1}{1 \vee (s \log(s))},
\hspace{ 1,3 cm}Q(k) \, := \, \int_0^k q(s) ds, \\
r(s,k) \, &:= \, \int_{(s \vee 0) \wedge k} ^k \frac{q(s')}{Q(k)} ds', \qquad  
\tn(s) \, := \,  \tau \, r(s-R,n-R).
\end{split}
\end{displaymath}
\begin{figure}[!htb] 
\begin{center}
\psfrag{1}{$0$}
\psfrag{2}{$R$}
\psfrag{3}{$n$}
\psfrag{4}{$s$}
\psfrag{5}{$\tau$}
\psfrag{6}{$\tn(s)$}
\includegraphics[scale=0.4]{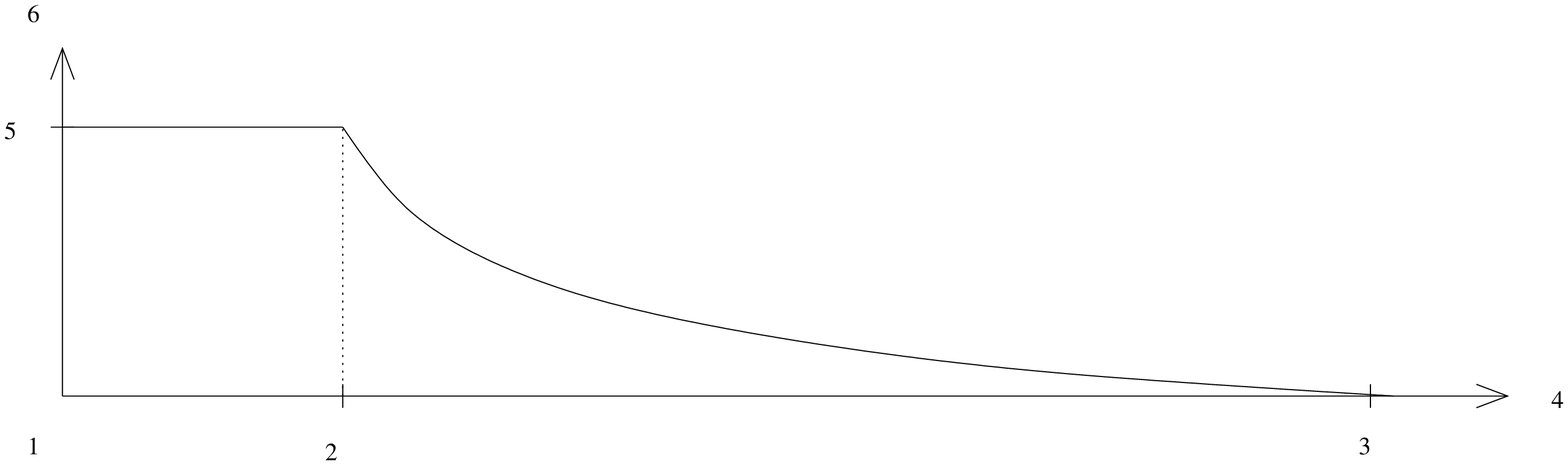}
\end{center}
\caption{Graph of $\tn$}\label{figtauhc}
\end{figure}
For a sketch of the graph of $\tn$ see Figure \ref{figtauhc}. 
Some properties of $\tn$ are:
\begin{equation} \label{tauinnenaussen}
\begin{split}
\tn(s) = \tau \, &\text{ for }  s \le R,  \quad 
\tn(s) = 0 \, \text{ for }  s \ge n, \quad  \text{ $\tn$ is decreasing.}
\end{split}
\end{equation}
Now $T_{R,n}(y) := y + \tn (|y|) \einh$
defines a transformation on $\RS$. This transformation  
can also be viewed as a transformation on $\Y$, 
such that every point $y$ of a configuration $Y$ is translated the distance 
$\tn(|y|)$ in direction $\einh$. We would like to use $T_{R,n}$ as a tool 
for our proof just as in \cite{FP1} and \cite{FP2}, but 
in order to deal with the hard core and the perturbation $u$, 
which is encoded in the bond process, we have to allow 
the transformation of a particle $y$ to depend on the 
configuration of particles in the neighborhood of $y$ and
the configuration of bonds joining $y$ to other particles. 
We thus aim to construct a transformation  
\[
\Tn: \Y \times \E \to \Y \times \E
\]
that is required to have the following properties: 
\begin{enumerate}
\item [(1)] For $B \subset E(Y)$ the transformed configuration  
$(\tY,\tB) = \Tn(Y,B)$  is constructed by translating every particle 
$y \in Y$ by a certain distance in direction  $\einh$, 
and by translating bonds along with the corresponding particles. 
\item [(2)] Particles in the inner region $\La_{n'-1}$ are translated by $\tau \einh$,
and particles in the outer region  ${\Lan}^c$ are not translated at all.
\item [(3)] Particles connected by a bond in  $B$ are translated the same 
distance.  
\item [(4)] $\Tn$ is bijective, and the density of the transformed process with 
respect to the untransformed process under the measure $\nu \otimes \pinh$ 
can be calculated. 
\item [(5)] We have suitable estimates on this density and on 
$H^{\bU}_{\Lan}(\tY)-H^{\bU}_{\Lan}(Y)$. For the last assumption we need particles 
within hard core distance to remain within hard core distance and 
particles at larger distance to remain at larger distance. 
\end{enumerate} 
Property (2) implies that the translation of the chosen cylinder event $D$ 
is the same as the transformation of $D$ by $\Tn$. Properties (3)-(5) are
chosen with a view to the right hand side of \eqref{desint}: If $\Tn$ has
these properties then the density of the transformed process with respect to the untransformed process under the measure $\mu \otimes \pin$ can be estimated. Therefore a transformation with these properties seems to be the right tool for proving \eqref{georgiianalog}. However, in general it is impossible to construct a transformation with all the given properties.  For example properties (2) and (5) cannot both be satisfied if $Y$ is a configuration of densely packed 
hard-core particles, or  properties (2) and (3) cannot both be satisfied 
if the inner and the outer region are connected by bonds. Similar problems arise 
for some of the other properties, so we will content ourselves 
with a transformation satisfying the above properties only for  
configurations $(Y,B)$ from a set of good configurations $\Gn$, which will
be shown to have probability close to $1$ for suitably chosen $R$ and $n$
in  Lemma~\ref{lebadsmall}. We define $\Gn$ to be of the form
\begin{equation} \label{good}
 \Gn  \, := \,\big\{ (Y,B) \in \Y \times \E: 
\; B \subset E(Y), \ryb <  R, \,
\sum \limits_{i=1}^5 \Sigma_i <   \de \big\}, 
\end{equation} 
where $\de \in ]0,1/2[$ is the constant chosen in section \ref{constants}.
The functions $\Sigma_i = \Sigma_i(R,n,Y,B)$ will be defined whenever we want good 
configurations to have certain properties, see \eqref{sig1}, \eqref{sig4} and \eqref{sig2}.
The condition involving $\ryb$ is meant to ensure that both the particle density 
and the number of bonds is not too high. More precisely for $Y \in \Y$ and 
$B \subset E(Y)$ let 
\[
B_+ := B \cup \{ y_1y_2 \in E(Y):  (y_1,y_2) \in \Kopp \}
\]
be the enlargement of $B$ by additional bonds between particles that are close 
to each other. We then define 
\[
\ryb \, := \, \sup\{ |y'| : y' \in C_{Y,B_+}(\La_{n'})\} 
\]
to be the range of the $B_+$-cluster of the inner region $\La_{n'}$. 
For $y \in Y$ let
\[
\ayb(y) \, := \, \min \{ \tn(|y'|) : y' \in  C_{Y,B_+}(y) \}.
\]
As $Y_{\Lan}$ is finite and 
$\tn(|.|) = 0$ on $\Lan^c$ by \eqref{tauinnenaussen}, this minimum is attained. 
By definition 
\begin{equation} \label{rnan}
(Y,B) \in \Gn, \;  y \in Y_{\La_{n'}} \quad \Rightarrow \quad \ayb(y) = \tau. 
\end{equation}


\subsection{Modifying the generalized translation}

We now define a transformation $\Tn$ with the properties described 
in the last section. As $n > R \ge n'$ are fixed throughout this 
section, we usually will omit the dependence on $n$ and $R$ in our notations. 
With a view to properties (1) and (3) we define the transformation
\[
\begin{split}
&\Tn: \Y \times \E \to \Y \times \E, \quad \Tn(Y,B) \, := \, (\Tnb(Y),\Tny(B))\\   
\text{ by } \quad 
&\Tnb(Y) \, := \, \bigcup_{k=0}^{m(Y,B)} (\cyb_k + \tauyb_k \einh) =
\{ y + \tyb(y) \einh : y \in Y \} \\
\text{ and }  \quad  
&\Tny(B) \, := \, \{ (y + \tyb(y)\einh)(y' + \tyb(y')\einh): yy' \in B \} 
\end{split}
\]
if $B$ is a subset of $E(Y)$, and  $\Tnb = id$ and $\Tny = id$ otherwise. 
Here $\big(\cyb_k\big)_{0 \le k \le m(Y,B)}$ is a certain partition of $Y$, where every $\cyb_k$ is a union of $B$-clusters.
\begin{figure}[!htb] 
\begin{center}
\psfrag{1}{$\tauyb_0 = 0$}
\psfrag{2}{$\tauyb_0 = 0$}
\psfrag{3}{$\tauyb_1$}
\psfrag{4}{$\tauyb_2$}
\psfrag{5}{$\tauyb_3$}
\psfrag{6}{$\tauyb_4$}
\psfrag{7}{$\tauyb_5 = \tau$}
\psfrag{a}{$\Lan$}
\psfrag{b}{$\La_R$}
\includegraphics[scale=0.65]{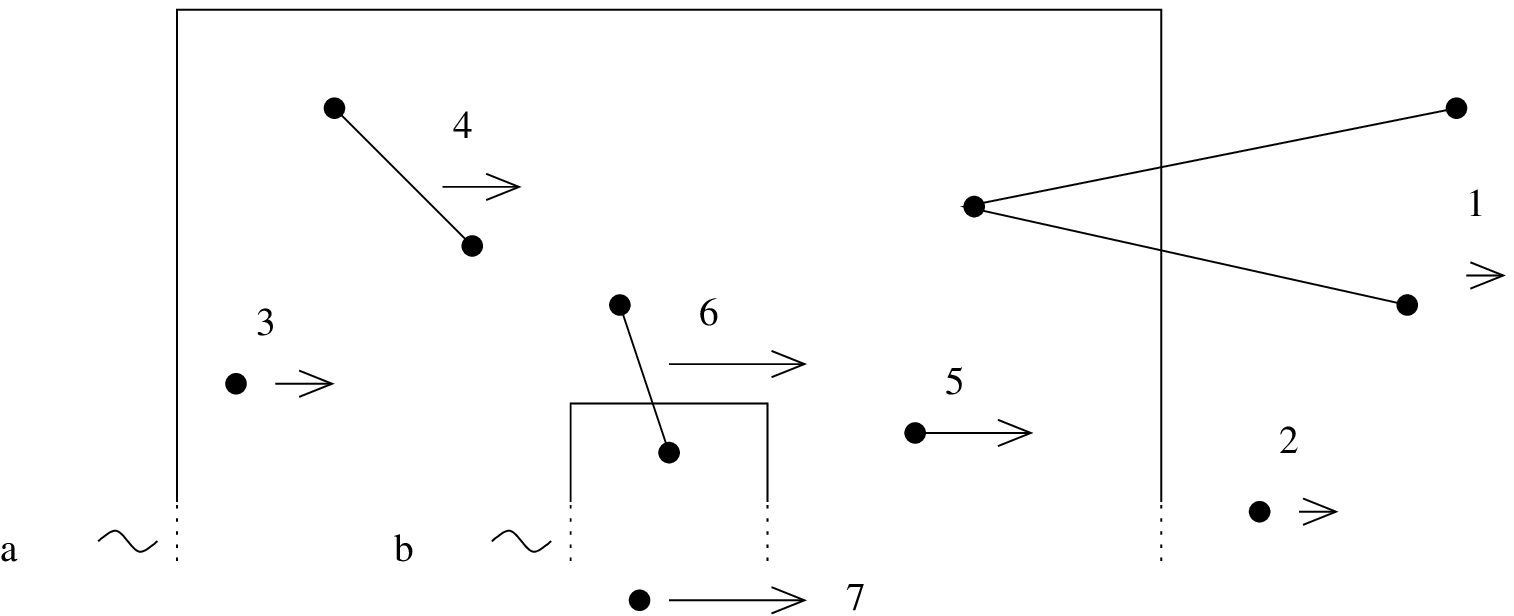}
\end{center} 
\caption{Every set $\cyb_k$ is translated by $\tauyb_k\einh$} 
\label{figclus}
\end{figure}
$\tauyb_k$ is the translation 
distance of all points in $\cyb_k$, and the translation distance function 
$\tyb: Y \to \R$ is defined by $\tyb(y) := \tauyb_k$ for $y \in \cyb_k$.
We are left to identify the points of $\cyb_k$ and their 
translation distances $\tauyb_k$. In our construction we would like 
to ensure that the sets $\cyb_k$ are ordered in a way such that 
\begin{equation} \label{mono}  
\tauyb_0 \le \tauyb_1 \le \ldots \le \tauyb_{m}. 
\end{equation}
This relation will be an important tool for showing the bijectivity of the 
transformation as required in property (4) of the last subsection. 
As required in (5) we also would like to have
\begin{align}
&y_1,y_2 \in Y,\,  (y_1,y_2) \in \Ko  \;
\Rightarrow \; \tyb(y_1) \,= \,\tyb(y_2), \label{vercon}\\
&y_1,y_2 \in Y, \, (y_1,y_2) \notin \Ko \; \Rightarrow \; 
(y_1+ \tyb(y_1)\einh,y_2+ \tyb(y_2)\einh) \notin \Ko. \label{greps} 
\end{align}
With these properties in mind we will now give a recursive definition 
of $\cyb_k$ and $\tauyb_k$ for a fixed $(Y,B) \in \Y \times \E$, 
where $B$ is a subset of $E(Y)$. 
In the $k^{th}$ construction step  $(k \ge 0)$ let 
\[
\begin{split} 
&\tyb_k \, := \,\tyb_{0} \wedge \bigwedge_{0 \le i < k} m_{\cyb_{i},\tauyb_{i}}\, 
= \, \tyb_{k-1} \wedge m_{\cyb_{k-1},\tauyb_{k-1}},\\
&\text{ where } \quad \tyb_0 \, := \, \tn(|.|) \quad  \text{ and } \quad 
m_{\cyb_i,\tauyb_i} \, := \,  \bigwedge_{y \in \cyb_i}  m_{y,\tauyb_i}. 
\end{split} \]
\begin{figure}[!htb]
\begin{center}
\psfrag{1}{$\La_R$}
\psfrag{2}{$\Lan$}
\psfrag{3}{$\epyb_2$}
\psfrag{4}{$\epyb_1$}
\psfrag{a}{$\epyb_3$}
\psfrag{6}{$\etauyb_1$}
\psfrag{5}{$\etauyb_2$}
\psfrag{b}{$\etauyb_3$}
\psfrag{8}{$m_{\ecyb_1,\etauyb_1}$}
\psfrag{7}{$m_{\ecyb_2,\etauyb_2}$}
\psfrag{c}{$m_{\ecyb_3,\etauyb_3}$}
\psfrag{9}{$\etyb_0$}
\psfrag{10}{$\tau$}
\psfrag{11}{$\R^2$}
\psfrag{12}{$\R$}
\includegraphics[scale=0.5]{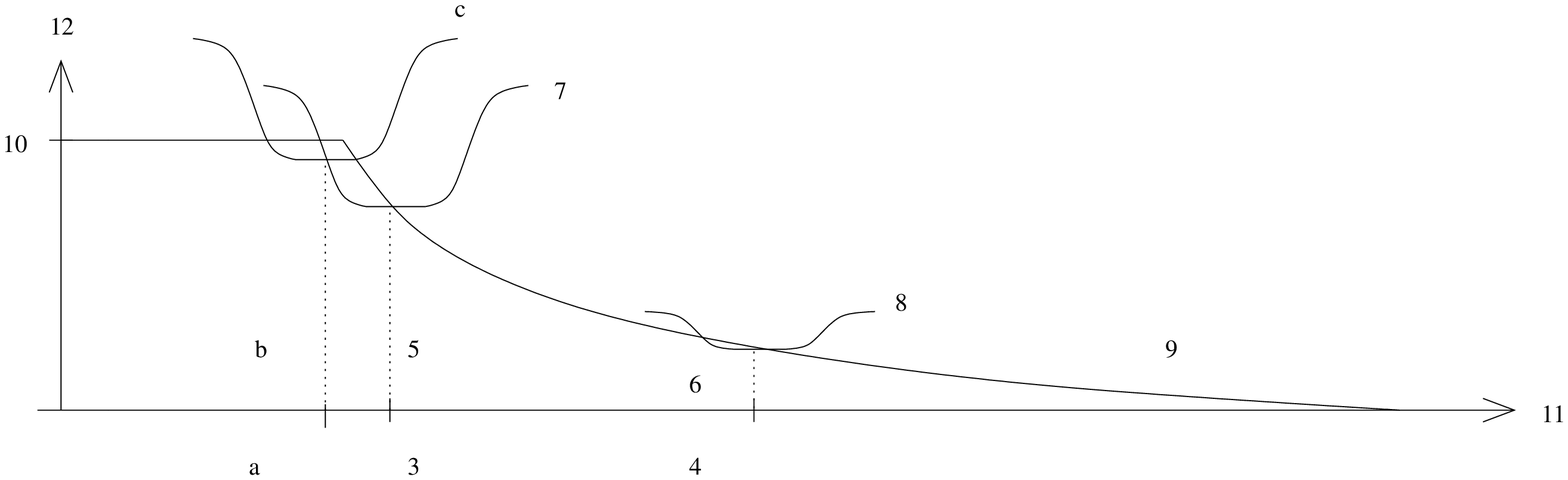}
\end{center}
\caption{Construction of $\etyb_k$ in the case that $\ecyb_k=\epyb_k$ is
one point only.}  \label{figtrdhc}
\end{figure}

\noindent
The auxiliary functions $\myt$ will be defined later. 
Let $\pyb_k$ be the set of points of 
$Y \weg (\cyb_0 \cup \ldots\cup \cyb_{k-1})$ 
at that $\tyb_k$ is minimal, and let $\tauyb_k$ be the corresponding minimal 
value, so that $\tauyb_k = \tyb_k(\pyb_k)$, i.e. $\tauyb_k = \tyb_k(y)$ 
for all $y \in \pyb_k$. Let $\cyb_k$ be the $B$-cluster of the set $\pyb_k$ and  $\Tyb_k := id + \tyb_k \einh$. The recursion stops when 
$Y \weg (\cyb_0 \cup  \ldots \cup \cyb_{m}) = \emptyset$, which occurs 
for a finite value of $m(Y,B) := m$. If it is clear from the context 
which configuration is considered, we may omit the dependence 
on $Y$ and $B$ in the above notations.\\
If $Y_{\Lan^c} \neq \emptyset$ then for every $y \in Y_{\Lan^c}$ we have 
$\tn(|y|) = 0$ by \eqref{tauinnenaussen}, so $y \in \epyb_0$ and 
$\etauyb_0 = 0$. 
This implies the second part of property (2). 
$\etyb_k$ is defined to be $\etyb_0 = \tn(|.|)$ modified by local distortions $\myt$. 
On the one hand we have thus ensured  that $\etyb_k-\etyb_0$ is small, 
i.e. $\etauyb_k \approx \tn(|y|)$ for all $y \in \epyb_k$, 
which will give us hold on the density in property (4). On the other hand 
the auxiliary functions of the form $\myt$ slow down the translation locally  
near every point $y'$ with known translation distance $t$, see Figure \ref{figtrdhc}. 
This will ensure properties \eqref{vercon} and \eqref{greps}. 
For $y' \in \RS$ and $t \in \R$ let the auxiliary function 
$\myt: \RS \to \baR$ be given by 
\[
\begin{split}
\myt(y) \, := \, \Bigg\{ 
&\begin{array}{cl}
t & \text{ if } \,  
 \hyt \cf > \frac{1}{2}\\
 t +  \hyt \fc(y',y)  + \infty \, 1_{\{ \fc(y',y)=1\}} \quad  
& \text{ otherwise}, 
\end{array}\\  
\text{ where } \quad &\hyt \, :=  \, |\tn(|y'| - \cK) - t|.
\end{split} 
\] 
\begin{figure}[!htb] 
\begin{center}
\psfrag{1}{$\Ko(y',r_2,\si)$}
\psfrag{2}{$\Kopp(y',r_2,\si)$}
\psfrag{3}{$t$}
\psfrag{4}{$\infty$}
\psfrag{5}{$h_{y',t}$}
\psfrag{6}{$r_1 \in \R$}
\psfrag{7}{$\myt(r_1,r_2,\si)$}
\includegraphics[scale=0.35]{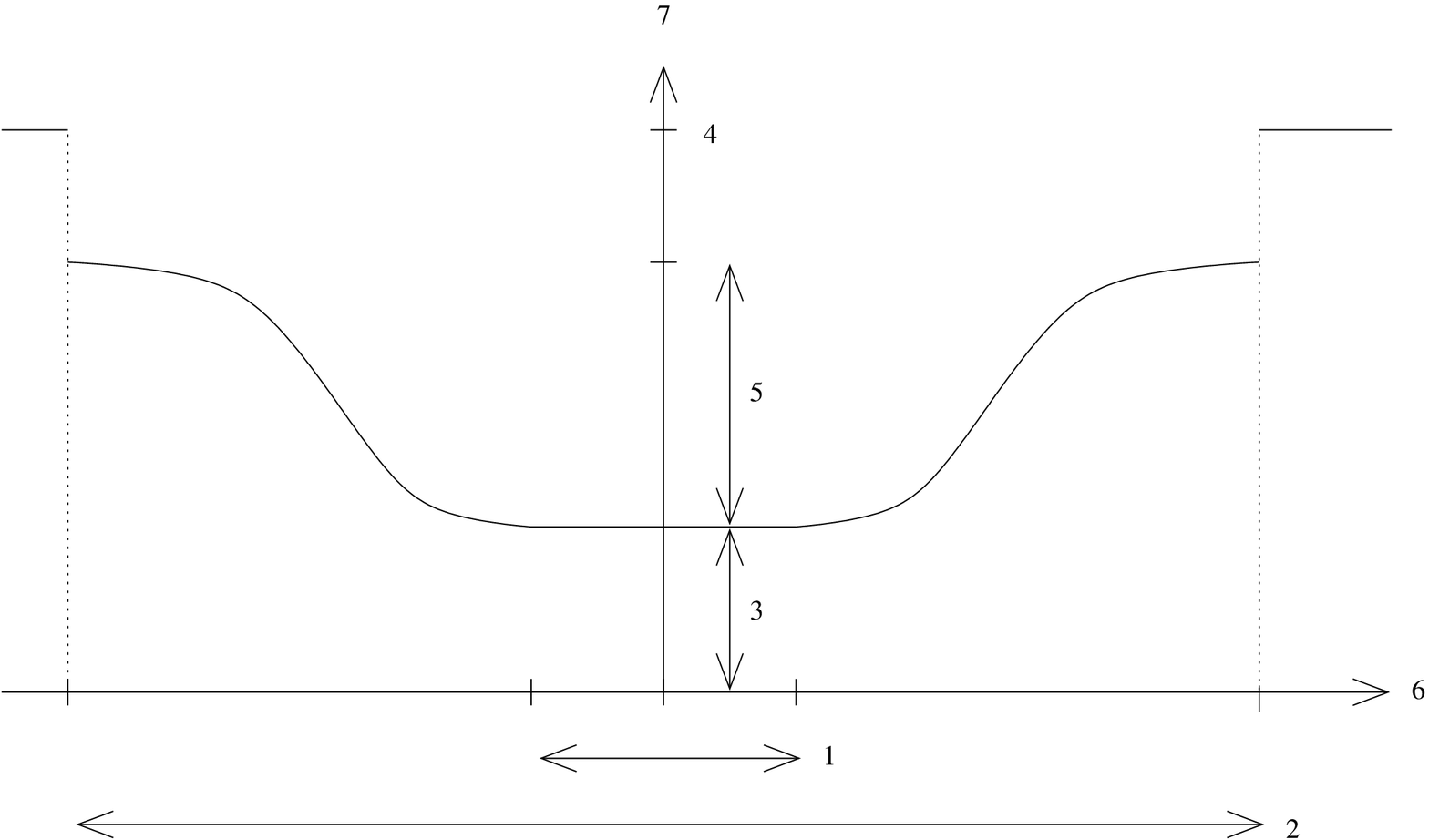}
\end{center}
\caption{Graph of $\myt(.,r_2,\si)$} \label{figmxthc}
\end{figure}
Note that the first case in the definition of $\myt$ has been introduced
in order to bound the slope of $\myt$. In Section~\ref{secmeig} we 
show important properties of this auxiliary function, but for the moment 
we will content ourselves with the intuition given by Figure \ref{figmxthc}.
Using Lemma~\ref{lemeas} one can show that all above objects are measurable 
with respect to the considered $\si$-algebras. In the rest of this section 
we will show that the above construction has indeed all the
required properties. 
\begin{lem}\label{lecons}
The construction satisfies \eqref{mono}, \eqref{vercon} and \eqref{greps}.
\end{lem}
\begin{lem} \label{leinnenaussen}
For good configurations  $(Y,B) \in \Gn$ we have    
\begin{equation} \label{innenaussen}
(\Tnb Y - \tau \einh)_{\La_{n'-1}} \, = \, Y_{\La_{n'-1}}\quad \text{ and }
\quad (\Tnb Y)_{{\Lan}^c} \, = \, Y_{{\Lan}^c}.
\end{equation}
\end{lem}
\begin{lem} \label{lebij}
The transformation $\Tn: \Y \times \E \to \Y \times \E$ is bijective.
\end{lem}
Actually in the proof of Lemma~\ref{lebij} we construct the inverse 
of $\Tn$. This is needed in the proof of Lemma~\ref{ledensp}, where 
we will show  for  every $\bY \in \Y$  that
$\nu_{\Lan} \otimes \pinh (.|\bar{Y})$ is absolutely continuous with respect to 
$\nu_{\Lan} \otimes \pinh (.|\bar{Y}) \circ \Tn^{-1}$ with density 
$\ph \circ \Tn^{-1}$, where 
\begin{equation} \label{densdef}
\ph (Y,B) \, := \, 
\prod_{k = 0}^{m(Y,B)} \prod_{y \in \pyb_k} \big|1 +  \partial_{\einh} \tyb_k (y)\big|. 
\end{equation}
Here $\partial_{\einh}$ is the spatial derivative in direction $\einh$.
The proof will also show that definition \eqref{densdef} 
makes sense  $\nu_{\Lan} \otimes \pinh (\,. \, |\bar{Y})$-a.s., in that 
all derivatives exist. We note that for every 
$y \in Y_{\Lan^c}$ the corresponding factor in the definition of $\ph(Y,B)$ equals $1$, so the above product in fact is a finite product. 
\begin{lem} \label{ledensp}
For every  $\bar{Y} \in \Y$ and every $\F_{\Y}  \otimes \F_{\E}$-measurable 
function  $f \ge 0$ 
\begin{equation} \label{dens}
\int d\nu_{\Lan}\otimes \pinh(.|\bar{Y}) \, (f \circ \Tn \cdot \ph) \,
= \, \int d\nu_{\Lan}\otimes \pinh(.|\bar{Y}) \, f.
\end{equation}
\end{lem}
Considering \eqref{georgiianalog} we also need  the backwards translation. So let
$\iTn$, $\iTnB$, $\iTny$, and $\bph$ be defined analogously to the above objects, 
where now $\einh$ is replaced by $-\einh$. The previous lemmas apply analogously 
to this deformed backwards translation. We note that $\iTn$ is not the inverse 
of $\Tn$.


\subsection{Final steps of the proof}  \label{Pott}

Let us now consider the $\delta >0$ and the test set $D$ chosen in 
Section~\ref{constants}. We identify $D$ with $D \times \E$ and 
use the shorthand notation $\Drn := D \cap \Gn$. With a view to Lemma~\ref{lekrit} we aim at 
showing that
\begin{equation} \label{geononne}
\mu \otimes \pin (\Tn \Drn) + \mu \otimes \pin (\iTn \Drn)  
 -  2(1-\de)\mu \otimes \pin(\Drn)
\end{equation}
is nonnegative. By \eqref{desint} and Lemma~\ref{ledensp} 
the first term of \eqref{geononne} equals 
\[ 
\begin{split}
&\int \mu(d\bY) \frac 1 {Z_{\Lan}^{U,z}(\bar{Y})} 
  \int \nu_{\Lan} \otimes \pinh(dY,dB|\bY) \;  1_{\Tn \Drn} \circ \Tn(Y,B)\\ 
&\hspace{1.5 cm} z^{ \# (\Tnb Y)_{\Lan}} \;  
  \ph(Y,B) \; e^{-H^{\bar{U}}_{\Lan}(\Tnb Y)} \prod_{b \in \Tny B} (e^{u(b)}-1).
\end{split} \]
By  Lemma~\ref{lebij} $~\Tn$ is bijective, 
by \eqref{innenaussen} we have $ \# (\Tnb Y)_{\Lan} =  \# Y_{\Lan}$ and 
by construction any two particles connected by a bond are translated the 
same distance. Hence the above integrand simplifies to 
\[
 1_{\Drn} (Y,B)  \, z^{\# Y_{\Lan}} \, 
  e^{\log \ph (Y,B) - H^{\bU}_{\Lan}(\Tnb Y)} \prod_{b \in B} (e^{u(b)}-1).
\]
Treating the other terms analogously, \eqref{geononne} can be seen to equal 
\[ \begin{split}
&\int \mu(d\bY) \frac{1}{Z_{\La_{n}}^{U,z}(\bar{Y})} 
  \int \nu_{\Lan} \otimes \pinh(dY,dB|\bY) \, 1_{\Drn} (Y,B)  
 \, z^{\# Y_{\Lan}}  \prod_{b \in B} (e^{u(b)}-1) \\  
& \times  \Big[
  e^{\log \ph(Y,B) -H^{\bU}_{\La_{n}}(\Tnb Y)} 
 + e^{\log \bph(Y,B) -H^{\bU}_{\Lan}(\iTnB Y)}  
 - 2(1-\de)e^{ -H^{\bU}_{\Lan}(Y)} \Big].
\end{split} \]
The convexity of the exponential function implies that the sum of the 
first two terms in the last bracket is greater or equal to 
\[
2 \, e^{\frac{1}{2}( \log \ph(Y,B) + \log \bph(Y,B)  
 -  H^{\bU}_{\La_{n}}(\Tnb Y) -  H^{\bU}_{\Lan}(\iTnB Y))},
\]
and here we can estimate the exponent using the following lemma: 
\begin{lem} \label{ngross}
For $(Y,B) \in \Gn$ we have     
\begin{equation} \label{dichtenklein}
\log \ph(Y,B)  + \log \bph(Y,B) \, \ge \, - \de \quad \text{ and } 
\end{equation}
\begin{equation} \label{taylor}
H^{\bU}_{\La_{n}}(\Tnb Y) + H^{\bU}_{\La_{n}}(\iTnB Y) \,
\le \, 2 H^{\bU}_{\Lan}(Y) + \de. 
 \end{equation}
\end{lem}
Using $e^{-\delta} \ge 1 - \delta$, this establishes the nonnegativity of
the above bracket and thus of \eqref{geononne}. So we have shown 
\begin{equation} \label{hingeorgii}
\mu \otimes \pin (\Tn \Drn) \, + \, 
\mu \otimes \pin (\iTn \Drn) \,
\ge \, 2\mu \otimes \pin (\Drn) -2\de.
\end{equation}
We would like to replace $\Drn$ by $D$. Using
 $D \in \F_{\Y,\La_{n'-1}}$ and \eqref{innenaussen} we obtain
\[
\alle (Y,B) \in \Drn: \quad (\Tnb Y - \tau \einh)_{\La_{n'-1}} \in D, 
\quad \text{ i.e. } \Tnb Y \in D + \tau \einh,
\]
and an analogous result for the backwards transformation. Hence 
\begin{equation} \label{traD}
\Tn(\Drn) \subset D + \tau \einh  \quad \text{ and } \quad 
\iTn(\Drn) \subset D - \tau \einh.
\end{equation} 
\begin{lem} \label{lebadsmall}
If the integers $n > R$ are chosen big enough, then $\mu \otimes \pin (\Gn^c) \le \de$.
\end{lem}
For the proof of Theorem~\ref{sym} we choose such $n > R$. 
Using \eqref{traD} and Lemma~\ref{lebadsmall} we deduce 
$\mu(D + \tau \einh) \, + \, 
\mu(D - \tau \einh) \, \ge \, 2\mu(D) \, - \,  4\de$
from \eqref{hingeorgii}. 
Taking the limit $\delta \to 0$, the claim of the theorem follows from Lemma~\ref{lekrit}.


\section{Proof of the lemmas from Section~\ref{secproofsym}} 
\label{secleproofsym}

\subsection{Convergence of energy sums: Lemma~\ref{umord}}

Let $n \in \N$. For every $Y \in \Y$ we have 
\begin{displaymath}
H^{\tiu}_{\Lan}(Y) \, = \, \sum_{b \in \en(Y)} \tiu(b) \, 
\le  \, \sideset{}{^{\neq}}\sum_{y_1, y_2 \in Y_{\Lan}} 
 \tiu(y_1,y_2) + 
\, \sum_{y_1  \in Y_{\Lan}} \sum_{y_2 \in  Y_{\Lan^c}} 
 \tiu(y_1,y_2), 
\end{displaymath}
and integrating this and applying Lemma~\ref{lekorab} for $\nu_{\Lan}(.|\bY)$ and $\mu$ we obtain 
\begin{displaymath}
\begin{split}
\int \mu \otimes \nu_{\Lan}(dY) &H^{\tiu}_{\Lan}(Y)  \, 
\le  \, \int_{\Lan} dy_1 \Big( \int_{\Lan} dy_2  \tiu(y_1,y_2)
+  z \xi \int_{\Lan^{\;c}} dy_2  \tiu(y_1,y_2) \Big)\\
&\le \, \int_{\Lan} dy_1 (1 + z\xi) \ex  \quad 
\le \quad 4n^2  (1 + z\xi) \ex \, \quad < \quad \infty, 
\end{split}
\end{displaymath}
where we have estimated the integrals over  $y_2$ by $\ex$ using  
\eqref{dec}. Thus we have proved the first 
assertion. However, $\mu$ is absolutely continuous with respect to 
 $\mu \otimes \nu_{\Lan}$, which follows from  \eqref{gibbsaequ} and the 
definition of the conditional Gibbs distribution. Hence the first 
assertion implies the second one.


\subsection{Properties of the auxiliary function } \label{secmeig}

A function $t: I \to \R$ on an interval $I$ is called \emph{hLd}, i.e. 
1/2-Lipschitz-continuous and differentiable at all but at most 
countably many points, if $|t(r)-t(r')| \le \frac 1 2  |r-r'|$ 
for all $r,r' \in I$, 
and if there is a countable set $M \subset I$ such that $f$ is 
differentiable in every point of $M \weg I$. The following lemmas 
show why we consider this type of function: 
\begin{lem} \label{hld}
Let $t:\R \to \R$ be hLd. Then the  transformation 
$T: \R \to \R$, $T := id + t$,  is bijective, strictly increasing, continuous, and 
differentiable a.e., and the Lebesgue transformation formula holds:
\begin{equation} \label{lebtra}
\int g(T(r)) T'(r) dr \, = \, \int g(r') dr' \quad \text{ for all 
measurable } g \ge 0. 
\end{equation}
\end{lem}

\Bew
We only need the 1/2-Lipschitz-continuity of $t$, which implies
\[
\frac 1 2 (r-r') \le T(r)-T(r') \, \le \, \frac 3 2  (r-r') \quad 
\text{ for all  } r \ge r' \in I,
\]
so $T$ is bijective, strictly increasing, and Lipschitz-continuous. 
The inverse $T^{-1}$ also is continuous and bijective, thus 
$\tla := \la \circ T$ is a measure on $(\R,\B)$. By the 
Lebesgue-Vitali differentiation theorem the Lipschitz-continuity of 
$T$ implies that $T$ is differentiable a.e. and 
$\frac{d\tla}{d\la} = T'$. Thus the transformation theorem 
implies \eqref{lebtra}. \qed

\begin{lem} \label{leanalmin} 
If $t_1,t_2: I \to \R$ are hLd functions on an interval $I$, 
then so is $t := t_1 \wedge t_2$, 
and we have $t'(s) \in \{t'_1(s) ,t'_2(s)\}$ whenever 
$t'(s)$ exists.  
\end{lem}

\Bew
The 1/2-Lipschitz-continuity of $t$ follows from the inequality
\[
\alle a_i,b_i \in \R: \quad |a_1 \wedge a_2 - b_1 \wedge b_2| 
\le |a_1-b_1| \vee |a_2-b_2|.
\]
For the differentiability let $M_i \subset I$ be a countable set 
such that $t_i$ is differentiable on $I \weg M_i$. Furthermore let 
\[
M_3 := \{r \in I \weg (M_1 \cup M_2): t_1(r) = t_2(r), 
t_1'(r) \neq t_2'(r) \}.
\] 
It is easy to check that every point of $M_3$ is isolated, so 
$M_3$ is countable. But $t_1 \wedge t_2$ is differentiable on 
$I \weg (M_1 \cup M_2 \cup M_3)$. Indeed, let 
$r \in I \weg (M_1 \cup M_2)$. If $t_1(r) \neq t_2(r)$, 
then $t$ coincides with one of the two functions in a neighborhood 
of $r$, and if $t_1(r) = t_2(r)$ and $t_1'(r) = t_2'(r)$, 
then $t$ is differentiable in $r$ with $t'(r) =t_1'(r) = t_2'(r)$. \qed\\

\noindent  
Let us call a function $t: \RS \to \R$  
$~\einh$-$1/2$-Lipschitz-continuous, $\einh$-differentiable,
or $\einh$-hLd if for all  $r_2 \in \R, \si \in S$ the function  $t(.,r_2,\si)$ has the corresponding property. 
\begin{lem} \label{lemeig}
For all $y' \in \RS$ and $t \in \R$ 
$~\tau_n(|.|) \wedge \myt$ is $\einh$-hLd.
\end{lem}   

\Bew
Let  $y' \in \RS$ and $t \in \R$. The claimed properties
concern the first spatial component only, so for fixed  $r_2 \in \R$ and 
$\si \in S$ we consider the functions
$\tilde{\tau} := \tau_n(|(.,r_2,\si)|)$, 
$\tilde{f} := \fc(y',.,r_2,\si)$, 
$\tilde{m}^f  :=   t +  \hyt \tilde{f}$, 
$\tilde{m} :=  \tilde{m}^f + \infty 1_{\{\tilde{f} =1\}}$.
It suffices to show that $\tilde{\tau} \wedge t$ is hLd for
$\hyt \cf > 1/2$, and $\tilde{\tau} \wedge \tilde{m}$ is hLd
for $\hyt \cf \le 1/2$. In order to get rid of the infinite part
of $\tilde{m}$ in the second case we define 
$I$ to be the convex hull of the closure of $\{ \tilde{f}<1 \}$. 
$I$ is a bounded closed interval, and we claim that  
\begin{equation} \label{taum}
\tilde{\tau} \wedge \tilde{m} = \tilde{\tau} \; 
\text{ on  } \; \overline{\R \weg I} \quad \text{  and } \quad 
\tilde{\tau} \wedge \tilde{m} = \tilde{\tau} \wedge \tilde{m}^f \; 
\text{ on  } \;I.
\end{equation}
Provided this is true, we are done by Lemma~\ref{leanalmin} as  $\tilde{m}^f$ is hLd for $\hyt \cf \le 1/2$ (by definition of $\cf$) and $\tilde{\tau}$ is hLd ($\tilde{\tau}$ is Lipschitz-continuous with Lipschitz-constant $\tau/Q(n-R) \le 1/2$, 
where we have used $n-R \ge 1$ and $\tau \le 1/2$). For a proof of 
\eqref{taum} we first observe that we have $\tilde{f} = 1$ on 
$\overline{\R \weg I}$ by the continuity of $\tilde{f}$ 
and thus $\tilde{m} = \infty$, which gives the first claim. 
For the second claim it suffices to show that for all $r \in I$ 
with  $\tilde{f}(r) = 1$ we have $\tilde{m}^f(r) \ge \tilde{\tau}(r)$. 
So let $r \in I$ with  $\tilde{f}(r) = 1$. 
We observe that $I$ is contained in the convex hull of the closure 
of $\Kopp(y',r_2,\si)$, as $\{ \tilde{f}<1 \} \subset \Kopp(y',r_2,\si)$
by definition of $\tilde{f}$. Thus for $y := (r,r_2,\si)$ 
we have $|y'-y| \le \cK$, which implies $|y'|-\cK \le |y|$. As 
$\tn$ is decreasing we obtain 
\[
\tilde{\tau}(r) \, = \, \tn(|y|) \, \le \, \tn (|y'|-\cK) \, 
\le \, t + \hyt \, = \, \tilde{m}^f (r) 
\]
by choice of  $\hyt$, and we are done.
 \qed


\subsection{Properties of the construction: Lemma~\ref{lecons}}

$\tyb_k$ is the minimum of finitely many functions of the form 
$\tau_n(|.|) \wedge \myt$, where $y' \in \RS$ and $t \in \R$. 
So the preceeding lemmas imply the following monotonicity and regularity properties of $\tyb_k$ and $\Tyb_k$:     
\begin{lem} \label{lelip}
For $Y \in \Y$, $B \subset E(Y)$ and $k \ge 0$ we have that  
$\tyb_k$ is $\einh$-hLd and  for every $y \in \RS$ $~\partial_{\einh} \tyb_k(y)$ equals $0$ or 
$\partial_{\einh} \tyb_0(y)$ or $\partial_{\einh} m_{y',\tyb(y')}(y)$
for some $y' \in Y$ with $(y,y') \in \Kopp$. Furthermore 
$\Tyb_k$ is $\lear$-increasing, $\einh$-continuous, 
bijective, and $\Tyb_k$ as a function of the first 
spatial coordinate satisfies \eqref{lebtra}.
\end{lem}
In the proofs of many of the following lemmas we need a way to 
calculate the translation distance of an arbitrary particle $y \in \cyb_k$
without knowing $\pyb_k$. This can be done using the following easy fact: 
\begin{lem} \label{lektk}
For $Y \in \Y$, $B \subset E(Y)$ and $k \ge 0$ we have 
\begin{equation}
\tauyb_k \, = \, \tyb_{k+1}(y) \quad \text{ for all }  y \in \cyb_k. \label{ktk}
\end{equation}
\end{lem}

\Bew 
For $y \in \ecyb_k$ we have $\etyb_{k+1}(y) 
= \etyb_k(y) \wedge \bigwedge_{ y' \in \ecyb_k} m_{y',\etauyb_k}(y)  
=  \etauyb_k$. Here we have used  $\etyb_k(y) \ge \etauyb_k$, which holds
by definition of  $\etauyb_k$,  $m_{y',\etauyb_k}(y) \ge \etauyb_k$, 
and $m_{y,\etauyb_k}(y) = \etauyb_k$. \qed\\

\noindent
For \eqref{mono} it suffices to observe that for every  $1 \le k \le m$ 
and $y \in \epyb_k$ we have 
\[
\etauyb_k \, = \, \etyb_k(y) \, 
= \, \etyb_{k-1}(y)  \wedge \bigwedge_{y' \in \ecyb_{k-1}} 
m_{y',\etauyb_{k-1}}(y) \, \ge \, \etauyb_{k-1}.
\]
This follows from the definition of $\etauyb_k$ and $\etyb_k$, 
from $\etyb_{k-1}(y) \ge \etauyb_{k-1}$ by the definition
of $\etauyb_{k-1}$, and from $\myt \ge t$. 
Now we will show \eqref{vercon} and 
\begin{equation} \label{grepss}
s \in [-1,1], y,y' \in Y, (y,y') \notin \Ko \; \Rightarrow \; 
(y,y'+ s(\tyb(y')-\tyb(y))\einh) \notin \Ko. 
\end{equation}
By the $\einh$-invariance of $\Ko$ \eqref{greps}
is equivalent to the special case $s=1$ of \eqref{grepss}.  
Let $y,y' \in Y$ and  $s \in [-1,1]$. Without loss of generality we may suppose 
that $y = y_i \in \ecyb_i$ and $y' = y_j \in \ecyb_j$, where  $0 \le i < j$.
(For $i=j$ we have $\tyb(y_i) = \tyb(y_j)$, so \eqref{vercon} and \eqref{grepss} are obvious.)  We now observe that 
$y_j \in \ela^i := \{y \in \RS: \etyb_i(y) \ge \etauyb_i\}$ and  
\begin{equation} \label{Kla}
\alle y \in \Ko(y_i) \cap \ela^i: \;\; 
\etyb_{j+1}(y) \, = \, \etyb_{i}(y) \wedge \bigwedge_{i \le k \le  j} \, 
  m_{\ecyb_k,\etauyb_k}(y) \, = \, \etauyb_i.
\end{equation}
This holds as $\etyb_i(y) \ge \etauyb_i$ by definition of $\ela^i$, 
$m_{\ecyb_k,\etauyb_k} \ge \etauyb_i$ by \eqref{mono}, and
$m_{\ecyb_i,\etauyb_i}(y) = \etauyb_i$ by $y \in \Ko(y_i)$.
If $(y_i,y_j) \in \Ko$, then $y_j \in  \Ko(y_i) \cap \ela^i$, 
so \eqref{Kla} and \eqref{ktk} imply   
$\etauyb_j = \etyb_{j+1}(y_j) = \etauyb_i$, which shows \eqref{vercon}.
For \eqref{grepss} suppose $(y_i,y_j) \notin \Ko$ and let 
$T_{j+1}^{s} := id + s \cdot \etyb_{j+1} \einh$. We have 
$y_j \in  \ela^i \weg \Ko(y_i)$ and $\etauyb_j = \etyb_{j+1}(y_j)$ by \eqref{ktk}, 
so it suffices to show that 
\begin{equation}  \label{grepst} 
T_{j+1}^{s}(\ela^i \weg \Ko(y_i)) \,
= \, \ela^i  \weg \Ko(y_i)  + s \etauyb_i \einh,  
\end{equation}
as this implies $y_j+s\etauyb_j \einh \notin \Ko(y_i)  + s \etauyb_i \einh$.
In order to show \eqref{grepst} we fix $\si \in S$ and $r \in \R$. Continuity of 
$\etyb_{i}(.,r,\si)$ implies $\etyb_i(.,r,\si) = \etauyb_i$ on 
$\partial \ela^i(.,r,\si)$. Just as in the proof \eqref{Kla} it follows that 
$\etyb_{j+1}(.,r,\si) = \etauyb_i$ on $\partial \ela^i$. But  
$T_{j+1}^{s}(.,r,\si)$ is increasing, continuous, and bijective, 
which can be shown as in the proof of Lemma~\ref{lelip}. So 
$T_{j+1}^{s}(\ela^i) \, = \, \ela^i + s \etauyb_i \einh$,  
and combining this with  \eqref{Kla} we are done.


\subsection{Properties of the deformed translation: Lemma~\ref{leinnenaussen}}

The following lemma shows how to estimate translation distances of particles.
\begin{lem} \label{leeigkons}
For good configurations $(Y,B) \in \Gn$ we have
\begin{equation}  \label{abdeftx}
\alle y \in Y: \quad 
0 \, \le \, \ayb(y) \, \le \, \tyb(y) \, \le \, \tn(|y|) \, \le \, \tau.
\end{equation}
\end{lem}
\Bew
The first and fourth inequality are a consequence of 
\eqref{tauinnenaussen}, and for the third it suffices 
to observe that for $y \in \ecyb_k$ we have 
$\etauyb_k \le \etyb_k(y) \le  \etyb_0(y)$ by the definition of $\etauyb_k$. 
For the second inequality we define
\begin{align} 
&\Sigma_1(R,n,Y,B) \, := \, \sum_{y,y' \in Y} 1_{\{ y \stackrel{Y,B_+}{\longleftrightarrow} y'\}} \Taun(y,y') 4 \cf^2, \label{sig1} \\
\text{ where } \quad &\Taun(y,y') \, := \, 1_{\{|y|\le |y'|\}} |\tn(|y|-\cK) - \tn(|y'|)|^2.
\label{Taun}
\end{align}
We have  $\Sigma_1(R,n,Y,B) <  1$ by  $(Y,B) \in \Gn$ and by definition of the set $\Gn$ in \eqref{good}. Hence every summand of 
$\Sigma_1$ is $<1$, and if we choose $y'$  to be a particle 
in $C_{Y,B_+}(y)$ such that $\tn(|.|)$ is minimal and $|y'| \ge |y|$ 
this implies 
\begin{equation} \label{verzer}
\forall y \in Y: \; 2 \cf \big|\tn(|y| - \cK) - \ayb(y)\big|  
\le \, 1.
\end{equation}
We will use this to show that all distortion functions  
$\myt$ in the definition of $\tyb(y)$ only have local influence 
in that in the definition of $\myt$ we have the second case 
($\hyt \cf \le 1/2$), which is needed in the following 
proof of
\[
 \ayb(y) \, \le \, \tauyb_k \quad \text{for all } y \in \cyb_{k}
\]
by induction on $k$. For  $k = 0$ we have equality. 
For the inductive step $k-1 \to k$ let  $i \le k-1$.  
By the third inequality of \eqref{abdeftx}, the inductive hypothesis, 
and  \eqref{verzer} we observe that for all $y_i \in \ecyb_i$ we have 
\[
0 \, \le \, \big(\tn(|y_i|-\cK) - \etauyb_{i}\big) \cf  \,  \le \, 
\big(\tn(|y_i| - \cK) - \ayb(y_i)\big) \cf \, 
\le \, 1 / 2, 
\]
so $h_{y,\etauyb_i} \cf \le 1/2$. Therefore $m_{y_i,\etauyb_i}$ is local in that $m_{y_i,\etauyb_i}(y) = \infty$
for all $y \in \epyb_k$ such that $(y_i,y) \notin \Kopp$. Thus
\[ 
\etauyb_k \,=\, \etyb_k(y)\, =  \, \etyb_0(y) \wedge \!\!\!\!
\bigwedge_{ y_i \in \ecyb_i: i<k,(y_i,y) \in \Kopp} \!\!\! m_{y_i,\etauyb_i}(y)  \, 
\ge \, \ayb(y),
\]
where the last step follows from 
$m_{y_i,\etauyb_i}(y) \ge \etauyb_i \ge \ayb(y_i)$, which is 
due to  the induction hypothesis, and from
$y_i \stackrel{Y,B_+}{\longleftrightarrow} y$ for $(y_i,y) \in \Kopp$. 
\qed\\

\noindent
We note that the proof Lemma~\ref{leeigkons} also shows that 
in the construction of $\Tn(Y,B)$ for a good configuration $(Y,B) \in \Gn$
all appearing distortion functions $\myt$ only have local influence 
as $\hyt \cf \le 1/2$. 
Now we will prove  Lemma~\ref{leinnenaussen}. It suffices to show for all 
$(Y,B) \in \Gn$ and $y \in Y$ that  
\begin{equation} \label{inaus}
\begin{split}
&y \in \La_{n'} \; \Rightarrow  \; \tyb(y) = \tau, \qquad 
y \in \La_{n'}^c \; \Rightarrow \; y + \tyb(y) \einh - \tau \einh \notin \La_{n'-1}\\
&y \in \Lan^c  \; \Rightarrow \;  \tyb(y) = 0, \qquad \text{ and } \quad  
y \in \Lan  \; \Rightarrow \; y + \tyb(y)\einh \in \Lan. 
\end{split}
\end{equation}
So let $(Y,B) \in \Gn$ and $y \in Y$. 
The first assertion of \eqref{inaus} now follows from \eqref{rnan} 
and \eqref{abdeftx}. 
The second assertion is an immediate consequence of $0 \le \tau- \tyb(y) \le 1$, 
which follows from \eqref{abdeftx} and $\tau \le 1$. 
The third assertion follows from \eqref{abdeftx} and \eqref{tauinnenaussen}, 
and for the fourth assertion let  $y \in \Lan$. As 
\[
y  \lear  y + \tyb(y)\einh  \lear  \Tyb_0(y)
\]
by  \eqref{abdeftx}, it suffices to show that also $\Tyb_0(y) \in \Lan$. 
This however follows from   $\Tyb_0 = id$ on ${\Lan}^c$ and the bijectivity 
of  $\Tyb_0$ from Lemma~\ref{lelip}.


\subsection{Bijectivity of  the transformation: Lemma~\ref{lebij}}
\label{secbij}

\begin{figure}[!htb] 
\begin{center}
\psfrag{1}{$\epyb_0$}
\psfrag{2}{$\etpyb_0$}
\psfrag{5}{$\epyb_2$}
\psfrag{6}{$\etpyb_2$}
\psfrag{7}{$\epyb_1$}
\psfrag{8}{$\etpyb_1$}
\psfrag{a}{$\ecyb_0$}
\psfrag{b}{$\etcyb_0$}
\psfrag{e}{$\ecyb_2$}
\psfrag{f}{$\etcyb_2$}
\psfrag{g}{$\ecyb_1$}
\psfrag{h}{$\etcyb_1$}
\psfrag{A}{$\etauyb_0$}
\psfrag{B}{$\ettauyb_0$}
\psfrag{E}{$\etauyb_2$}
\psfrag{F}{$\ettauyb_2$}
\psfrag{G}{$\etauyb_1$}
\psfrag{H}{$\ettauyb_1$}
\psfrag{x}{$\Tn: (Y,B) \mapsto (\tY,\tB)$}
\psfrag{y}{$\tTn: (\tY,\tB) \mapsto (Y,B)$}
\includegraphics[scale=0.5]{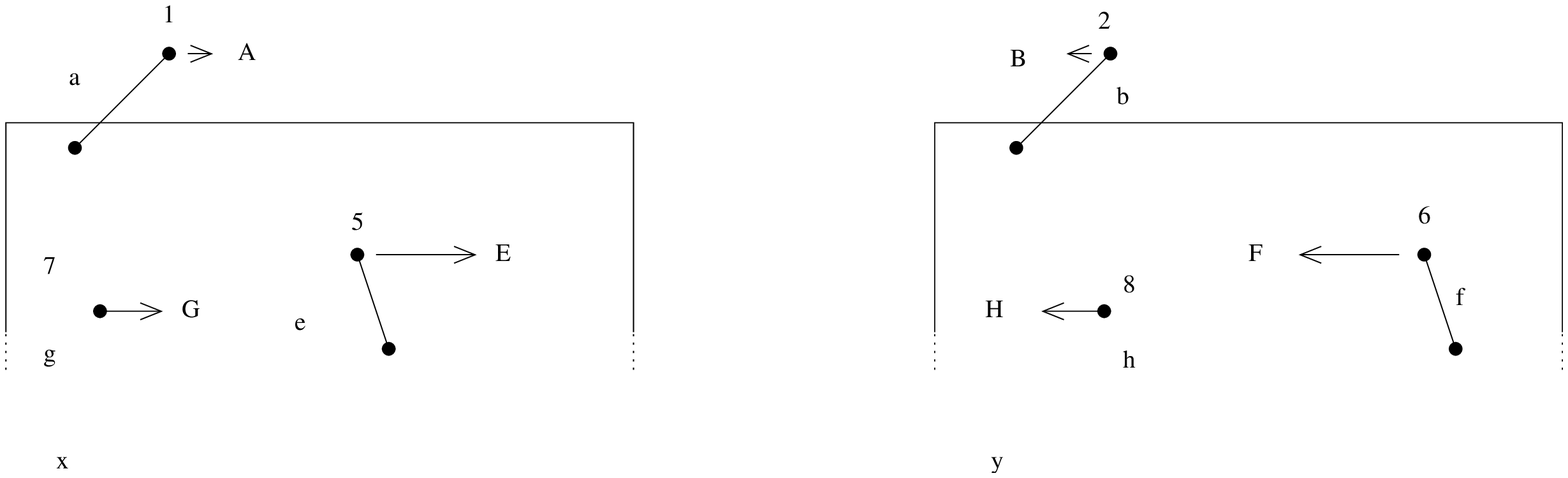}
\end{center}
\caption{Construction of the inverse  $\tTn$ of $\Tn$.} \label{figinv}
\end{figure}

We construct the inverse transformation $\tTn$ recursively, 
similarly to the construction of $\Tn$, 
i.e. from a given configuration $(\tY,\tB)$ we will choose 
sets of points $\tcyb_k$ and translate them by $\ttauyb_k$ 
in direction $-\einh$. To get an idea how to define the inverse
transformation we start with a fixed configuration $Y \in \Y$, 
$B \subset E(Y)$ and set $(\tY,\tB) := \Tn(Y,B)$. 
In the construction of $\Tn(Y,B)$ 
we defined a partition of $Y$ into sets of 
particles $\ecyb_k$, corresponding sets $\epyb_k$, and translation distances 
$\etauyb_k$. We denote the corresponding image sets by 
$\etpyb_k := \epyb_k + \etauyb_k \einh$ and $\etcyb_k := \ecyb_k + \etauyb_k \einh$, 
see Figure \ref{figinv}.
For the construction of the inverse transformation we have to find a method to
identify the sets $\etcyb_k$ among the points of $\tY$ without knowing $(Y,B)$. 
Suppose now that we have already found $\etcyb_i$, $\etpyb_i$ and 
$\etauyb_i$ for all $i<k$, then  with this information we are able to reconstruct $\ecyb_i := \etcyb_i - \etauyb_i \einh$ 
and thus $\etyb_k$ and $\eTyb_k$.
The following lemma tells us, how to find $\etpyb_k$ in that case, which 
will also give us $\etcyb_k$, $\epyb_k := \eTyb_k^{-1}(\etpyb_k)$
and $\etauyb_k := \etyb_k(\epyb_k)$. 
\begin{lem} \label{leminim}
Let $0 \le k \le m$. $\etpyb_k$ is the set of points of  
$\tY \weg \bigcup_{i \le k-1} \etcyb_{i}$ where 
$\etyb_k \circ \eTyb_k^{-1}$ attains its minimum value. 
\end{lem} 
\Bew 
We first observe that for all $k$  by definition of $\etyb_k$ we have
\begin{equation} \label{trdtinv}
\eTyb_k^{-1} + \etyb_k \circ \eTyb_k^{-1}\einh \, = \, id.  
\end{equation}
Since $\etyb_{k+1} \le \etyb_{k}$, we also have $\eTyb_{k+1} \lear \eTyb_{k}$, 
and therefore  $\eTyb_{k}^{-1} \lear  \eTyb_{k+1}^{-1}$ by the $\einh$-monotonicity  of $\eTyb_{k+1}^{-1}$ from Lemma~\ref{lelip}. 
Together with \eqref{trdtinv} this implies 
\begin{equation} \label{inveincr}
\etyb_{k+1} \circ \eTyb_{k+1}^{-1} \, \le \, \etyb_k \circ \eTyb_k^{-1}.
\end{equation}
Now let $0 \le k \le m$, $\ty_k \in \etpyb_k$, and  $\ty_l \in \etcyb_l$
for some $l \ge k$. Then $y_k := \ty_k - \etauyb_k \einh \in \epyb_k$
and  $y_l := \ty_l - \etauyb_l \einh \in \ecyb_l$. By definition and by 
\eqref{ktk} we have $\eTyb_{k}(y_k) \, = \, \ty_k$, 
$\etyb_{l+1}(y_l) \, = \, \etauyb_l$, and $\eTyb_{l+1}(y_l) \, = \, \ty_l$. 
Using \eqref{mono} and \eqref{inveincr} we deduce 
\[
 \etyb_{k} (\eTyb_{k}^{-1} (\ty_k))\, 
= \, \etauyb_k\, \le \, \etauyb_l \, 
= \, \etyb_{l+1} (\eTyb_{l+1}^{-1} (\ty_l)) \,
\le \, \etyb_{k} (\eTyb_{k}^{-1} (\ty_l)).
\]
If for the given $\ty_l$ we have equality, all inequalities 
in the previous line have to be equalities, so $\etauyb_k = \etauyb_l$  and 
$\etyb_{l+1} (\eTyb_{l+1}^{-1} (\ty_l)) =\etyb_{k} (\eTyb_{k}^{-1} (\ty_l))$. 
Combining this with \eqref{trdtinv} we get 
$y_l = \eTyb_{l+1}^{-1} (\ty_l) =\eTyb_{k}^{-1} (\ty_l)$, 
so $\eTyb_{k}(y_l)= \ty_l$ and thus $\etyb_k(y_l) = \etauyb_l = \etauyb_k$. 
By definition of $\epyb_k$ we conclude $y_l \in \epyb_k$, 
so $\ty_l \in \etpyb_k$ and we are done.  \qed \\

\noindent
Lemma~\ref{leminim} tells us exactly how to construct the inverse of $\Tn$ 
recursively. So let $\tY \in \Y$ and $\tB \subset E(\tY)$.
In the $k^{th}$ construction step $(k \ge 0)$ we define 
\[
\ttyb_{k} \, :=  \, \ttyb_{k-1} \wedge  \! \! \bigwedge_{y \in \tcyb_{k-1}-\ttauyb_{k-1}} 
  \! \! m_{y,\ttauyb_{k-1}}, \quad \text{ where } \quad \ttyb_0 := \tn(|.|).
\]
Let $\tTyb_k = id + \ttyb_k \einh$, and define $\tpyb_k$ 
to be the set of particles
of $\tY \weg (\tcyb_0 \cup \ldots  \cup \tcyb_{k-1})$ at that the minimum
of $\ttyb_k \circ (\tTyb_k)^{-1}$ is attained. Let  
$\ttauyb_k := \ttyb_k \circ (\tTyb_k)^{-1}(\tpyb_k)$ be  the corresponding 
minimal value and $\tcyb_k$ be the $\tB$-cluster of the set $\tpyb_k$. The recursion stops, when  
$\tY \weg (\tcyb_0 \cup \ldots  \cup \tcyb_{\tm}) = \emptyset$, 
which will occur for a finite value of $\tm(\tY,\tB) := \tm$. Again, sometimes we will omit the dependence on $\tY$ and $\tB$ in our notations if it is clear from the context which 
configuration is considered. 
We need to show that the above construction is well defined, i.e. 
that   $\tTyb_k$ is invertible in every step. Furthermore we 
need some more properties of the construction: 
\begin{lem} \label{leeiginv}
Let $\tY \in \Y$, $\tB \subset E(\tY)$ and $k \ge 0$. Then
\begin{align}
& \ttyb_k \text{ is $\einh$-hLd}, \quad 
\tTyb_k \text{ is bijective and $\lear$-increasing}, \label{lipinv}\\
& (\tTyb_k)^{-1} + \ttyb_k \circ (\tTyb_k)^{-1} \einh  \, = \,id, 
 \label{TSinv}\\
& \alle c \in \R ,y \in \RS: \, \ttyb_k \circ (\tTyb_k)^{-1}(y) 
 \,\ge\,c\, 
  \Leftrightarrow \, \ttyb_k(y -  c\einh) \, \ge \, c, \label{tsinv}\\ 
& \ttyb_{k} \, \le \, \ttyb_{k-1} \quad \text{ and } \quad  
 \ttauyb_{k-1} \, \le \, \ttauyb_{k}, \label{monoinv}\\
& \alle y \in \tcyb_{k}: \, \ttyb_{k+1} \circ(\tTyb_{k+1})^{-1}(y) \, 
= \, \ttauyb_k.
\label{kek}
\end{align} 
\end{lem}
\Bew
The definitions of $\ettyb_{k}$ and $\etTyb_k$ are similar to those of 
$\etyb_k$ and $\eTyb_k$, so we can show \eqref{lipinv} and \eqref{TSinv} 
just as the corresponding properties in Lemma~\ref{lelip}  and  
\eqref{trdtinv}. 
For \eqref{tsinv} we note that for $c \in \R$ and $y \in \RS$ 
the equivalence 
\[
\begin{split}
\ettyb_k \circ (\etTyb_k)^{-1}(y) \,\ge \,c \quad &\Leftrightarrow \quad  
(\etTyb_k)^{-1}(y)  \, \lear \, y - c\einh \\ 
&\Leftrightarrow \quad  
y \, \lear \, \etTyb_k(y - c\einh ) = y - c\einh + \ettyb_k(y - c\einh )
\einh
\end{split}
\]
follows from   \eqref{TSinv} and \eqref{lipinv}.
The first part of \eqref{monoinv} is obvious and for the second part  
we observe that for $\ty_k \in \etpyb_k$ we have 
\[
\begin{split}
&\ettyb_{k-1} \circ (\etTyb_{k-1})^{-1}(\ty_k) \ge \ettauyb_{k-1} 
\quad \Rightarrow \quad 
\ettyb_{k-1}(\ty_k - \ettauyb_{k-1}\einh ) \ge \ettauyb_{k-1}  \\ 
&\Rightarrow \quad 
\ettyb_{k}(\ty_k - \ettauyb_{k-1} \einh) \ge \ettauyb_{k-1}  
\quad \Rightarrow \quad 
\ettauyb_k =\ettyb_{k} \circ (\etTyb_k)^{-1}(\ty_k) \ge \ettauyb_{k-1},
\end{split} 
\]
where the first statement holds by definition of $\etpyb_{k-1}$, 
the first and the third implication hold by  \eqref{tsinv}, and the 
second implication is by definition of $\ettyb_k$.
For \eqref{kek} let  $\ty_k \in \etcyb_k$. We have
\[ \begin{split}
\ettyb_{k} &\circ \etTyb_{k}^{-1}(\ty_k) \,\ge \,\ettauyb_k 
\quad \Rightarrow \quad 
\ettyb_{k}(\ty_k-\ettauyb_k \einh ) \,\ge \,\ettauyb_k \\
&\Rightarrow \quad \ettyb_{k+1}(\ty_k - \ettauyb_k \einh) \,= \, \ettauyb_k 
\quad \Rightarrow \quad 
\ettyb_{k+1} \circ \etTyb_{k+1}^{-1}(\ty_k) \,= \, \ettauyb_k,
\end{split} \]
where the first statement holds by definition, and the implications 
follow from \eqref{tsinv}, 
$\ty_k - \ettauyb_k \einh \in \etcyb_k - \ettauyb_k \einh$ and 
\eqref{TSinv} respectively. \qed\\ 

\noindent
For every $0 \le k \le \tm(\tY,\tB)$ and $\ty_k \in \tcyb_k$ let 
$\ttyb(\ty_k) := \ttauyb_{k}$. This defines a translation distance function 
$\ttyb: \tY \to \R$. We define 
\[ 
\begin{split}
\tTnb(\tY) \, &:= \, \bigcup_{k=0}^{\tm(\tY,\tB)} (\tcyb_k - \ttauyb_k \einh) =
\{ y - \ttyb(y) \einh : y \in \tY \}  \quad \text{ and }\\ 
\tTny(\tB) \, &:= \, \{ (y - \ttyb(y)\einh)(y' - \ttyb(y')\einh): yy' \in \tB \}. 
\end{split}
\]
Now if $\tB$ is a not a subset of  $E(\tY)$ we define   
$\tTnb = id$ and $\tTny = id$. Let
\[
\tTn: \Y \times \E \to \Y \times \E, \quad \tTn(\tY,\tB) \, 
:= \, (\tTnb(\tY),\tTny(\tB)).  
\]
By Lemma~\ref{lemeas} we see again that all above objects are 
measurable with respect to the considered $\si$-algebras. 
The only difficulty is to show that $(\tTyb_k)^{-1}(y)$  
is measurable. This however follows from the  $\einh$-monotonicity of  
$\tTyb_k$.
In order to show that  $\tTn$ really is the inverse of $\Tn$ we need an 
analogue of Lemma~\ref{leminim}. Let  $\tY \in \Y$ and $\tB \subset E(\tY)$. 
Let $\ettyb_k$, $\etTyb_k$, $\etcyb_k$, $\etpyb_k$, and  $\ettauyb_k$ 
$~(0 \le k \le \tilde{m})$ as above and denote $(Y,B) := \tTn(\tY,\tB)$,
 $\epyb_k := \etpyb_k - \ettauyb_k \einh$, and 
$\ecyb_k := \etcyb_k - \ettauyb_k \einh$, see Figure~\ref{figinv}.
\begin{lem} \label{leminiminv}
Let $0 \le k \le \tm$. $\epyb_k$ is the set of points of 
$Y \weg \bigcup_{i \le k-1} \ecyb_i$ where  
$\ettyb_k$ attains its minimal value. 
\end{lem} 
\Bew 
Let $0 \le k \le \tm$, $y_k \in \epyb_k$, and $y_l \in \ecyb_l$ for some $l \ge k$. 
Then $\ty_l := y_l+ \ettauyb_l \einh  \in \etcyb_l$ and 
$\ty_k := y_k+ \ettauyb_k \einh  \in \etpyb_k$. By definition of 
$\ettauyb_k$ and $\ettauyb_l$, by  \eqref{kek} and \eqref{TSinv}
we have $\etTyb_{k}^{-1}(\ty_k) = y_k$, 
$\ettyb_{l+1} (\etTyb_{l+1}^{-1}(\ty_l))= \ettauyb_l$, and 
$\etTyb_{l+1}^{-1}(\ty_l) = y_l$.  
Thus from \eqref{monoinv} we deduce  
\[
\ettyb_{k} (y_k)\, 
= \,  \ettauyb_k\, \le \, \ettauyb_l \, 
= \, \ettyb_{l+1} (\etTyb_{l+1}^{-1} (\ty_l)) \,
= \, \ettyb_{l+1} (y_l) \, 
\le \, \ettyb_{k} (y_l).
\]
If for the given $y_l$ we have equality, all inequalities in the previous line
have to be equalities, so $\ettauyb_k = \ettauyb_l$ and 
$\ettyb_{k} (y_l) =\ettauyb_l$, i.e.   
$\etTyb_{k} (y_l) = y_l +  \ettauyb_l \einh = \ty_l$. 
This gives $\ettauyb_k = \ettauyb_l =  \ettyb_{k} (y_l) 
=  \ettyb_{k}(\etTyb_{k}^{-1} (\ty_l))$. By definition of 
$\etpyb_k$ we conclude $\ty_l \in \etpyb_k$, hence  $y_l \in \epyb_k$   
and we are done. \qed \\
\begin{lem}\label{leinv}
On $\Y \times \E$ we have $ \quad  \tTn \circ \Tn \, = \, id \quad $ 
and  $ \quad  \Tn \circ \tTn \, = \, id$.
\end{lem}
\Bew
For the first part let  $Y \in \Y$, $B \in \E$, and $(\tY,\tB) := \Tn(Y,B)$.
If $B$ is not a  subset of $E(Y)$  we have 
$\tTn \circ \Tn (Y,B) = \tTn  (Y,B)  = (Y,B)$ and we are done. 
Else it suffices to prove 
\begin{equation} \begin{split} \label{tildegleich}
\ttyb_k \, &= \, \tyb_k, \quad \tTyb_k \, = \, \Tyb_k, \quad 
\ttauyb_k \, = \, \tauyb_k, \\ 
\tpyb_k \, &= \, \pyb_k +  \tauyb_k \einh, \quad \text{ and } \quad  
 \tcyb_k \, = \, \cyb_k  +  \tauyb_k \einh
\end{split}\end{equation} 
for every $k \ge 0$ by induction on  $k$. For the inductive step $k-1 \to k$ 
we observe that $\ttyb_k = \tyb_k$ by induction hypothesis, 
and $\tTyb_k  =  \Tyb_k$ is an immediate consequence. Combining this 
with Lemma~\ref{leminim} and the definition of $\tpyb_k$ we get 
$\tpyb_k = \pyb_k  +  \tauyb_k \einh$ and   $\ttauyb_k = \tauyb_k$.
$\tcyb_k = \cyb_k +  \tauyb_k \einh$ is an immediate consequence.
The case $k=0$ can be shown similarly: Here  $\ttyb_0 = \tyb_0$ 
holds by definition and the rest again follows from Lemma~ \ref{leminim}.\\ 
For the second part let  $\tY \in \Y$, $\tB \in \E$, and 
$(Y,B) := \tTn(\tY,\tB)$. As above we may assume 
$\tB \subset E(\tY)$, and it suffices to show \eqref{tildegleich}
by induction on $k$. Here the inductive step follows from  Lemma~\ref{leminiminv}.
\qed


\subsection{Density of the transformed process: Lemma~\ref{ledensp}}

By definition the left hand side of \eqref{dens} equals 
\[ 
e^{-4 n^2} \sum_{k \ge 0} \frac{1}{k!} I(k), \quad
\text{ where }  \,
  I(k) \, = \, \int_{{\Lan}^{k}} dy
   \sideset{}{'}\sum_{B \subset \en(\bY_y)} (f \circ \Tn \cdot \ph) (\bY_y,B),
\]
where the summation symbol $\sum'$  indicates that the sum extends over finite subsets only, and we have used the shorthand notation 
$\bY_y = \{y_1,\ldots ,y_k\} \cup \bY_{\Lan^c}$
for $y \in (\Lan \times S)^k$.
We would like to fix the bond set $B$ before we choose the particle 
states $y_i$. Thus we introduce bonds between indices of 
particles instead of bonds between particles. Let  $\Nk := \{1,\ldots ,k\}$, 
$\bY^k := \Nk \cup \bY_{\Lan^c}$,  and 
$E_n(\bY^k) := \{ y_1y_2 \in E(\bY^k): y_1y_2 \cap \Nk \neq \emptyset \}$.
For $B \subset E_n(\bY^k)$ and $y  \in (\Lan \times S)^{I}$ 
$\,(I \subset \Nk)$ we define $B_y$ to be the bond set constructed from $B$ 
by replacing the point  $i \in I$ by $y_i$ in every bond of $B$ and by deleting 
every bond $B$ that contains a point $i \in \Nk \weg I$. Analogously let 
$\bY_y := \{ y_i: i \in I\} \cup \bY_{\Lan^c}$ be the configuration corresponding
to the sequence and let $(\bY,B)_y := (\bY_y,B_y)$. We obtain
\[
I(k)  \,=  \! \! \sideset{}{'}\sum_{B \subset E_n(\bY^k)} \! \! \! \! I(k,B),   
\; \text{ where } \,   I(k,B)  \,:= \,  \int_{{\Lan}^{k}} \! dy \,
   (f \circ \Tn \cdot \ph)(\bY,B)_y.   
\]
To compute $I(k,B)$ we need to calculate  $\Tn(\bY,B)_y$, and for this 
we must identify the points of $P_i^{\bY_y,B_y}$ among $\bY_y$.
So let $\Pi_k$ the set of all sequences $\per = (\per_j)_{0 \le j \le m}$
of disjoint nonempty subsets of $\bY^k$ such that $ \bY_{\Lan^c} \subset \per_0$
and every $B$-cluster of $(\bY^k,B)$ has nonempty intersection with exactly 
one of the sets $\per_j$, i.e. the $B$-clusters $\per_j^{B}$ of the sets $\per_j$
define a partition of $\bY^k$.  
Let the length of the sequence be denoted by $m(\per) := m$. For $\per \in \Pi_k$
and  $y  \in (\Lan \times S)^k$ let $\per_y =  (\per_{j,y})_{0 \le j \le m(\per)}$ and $\per^{B}_y=  (\per^{B}_{j,y})_{0 \le j \le m(\per)}$
be the sequences corresponding to $\per$ and $\per^B$, where every $i$ is replaced by $y_i$.  
For $\per \in \Pi_k$ let 
\[ \begin{split}
\akb  &:=  \big\{ y \in (\Lan \times S)^{k}: m(\per) = m(\bY_y,B_y), 
\alle j \ge 0: \per_{j,y} = P_j^{\bY_y,B_y}\big\}, \\
\takb  &:=  \big\{ y \in (\Lan \times S)^{k}: m(\per) = \tm(\bY_y,B_y), 
\alle j \ge 0: \per_{j,y} = \tilde{P}_j^{\bY_y,B_y}\big\},
\end{split}\]
where $\tm(\bY_y,B_y)$ and $\tilde{P}_j^{\bY_y,B_y}$ are the
objects  from the construction of the inverse transformation 
in Subsection~\ref{secbij}. We note that 
\[
\alle y \in \akb: \per^{B}_{j,y} = C_j^{\bY_y,B_y} \quad \text{ and } \quad 
\alle y \in \takb: \per^{B}_{j,y} = \tilde{C}_j^{\bY_y,B_y}.
\] 
Now we can write
\[
I(k,B) \, = \,\sum_{\per \in \Pi_k} \int dy \, 
        1_{\akb}(y) (f \circ \Tn \cdot \ph)(\bY,B)_y,   
\]
and we denote the summands in the last term by  $I(k,B,\per)$. 
If  $y \in \akb$ we can derive a simple expression for $\Tn(\bY,B)_y$: 
For $i \in \per^B_j$ the translation distance of $y_i$ doesn't 
depend on all components of $y$, but only on those $y_l$ such that 
$ l \in \per^B_{j'}$ for some $j' \le j-1$ and  additionally on  
those $y_l$ such that $ l \in \per_{j}$ whenever 
$i \notin \per_j$. Hence for  $y \in (\Lan \times S)^k$, $\per \in \Pi_k$ and 
 $0 \le j \le m(\per)$ we define $y^{\per,j}$ to be the subsequence of $y$ 
corresponding to the index set $\bigcup_{j' \le j} \per^{B}_{j'}$, we define 
a formal translation distance and a formal transformation by 
\[ 
\begin{split}
&t_j^{B,\per,y} := t^{(\bY,B)_{y^{\per,j-1}}}_{j}
\quad \text{ and } \quad  
T^{B,\per}(y) := (T_{j(i)}^{B,\per,y} (y_i))_{1 \le i \le k}, \quad \text{ where }\\ 
& j(i) := j \; \text{ for }  \;i \in \per_j \quad \text{ and } \quad 
 j(i) := (j,*)  \;\text{ for }  \;i \in \per^B_j \weg \per_j, \\
&T_j^{B,\per,y} := id + t_j^{B,\per,y} \; \einh \; \text{ and } \;  
T_{j,*}^{B,\per,y} :=  id + t_j^{B,\per,y}(y_{\min \per_j})\einh. 
\end{split}
\]
Then 
\begin{equation} \label{at}
y \in \akb \, \Rightarrow  \, \left\{ 
\begin{aligned} &\T_n(\bY,B)_y \, 
= \, (\bY, B)_{T^{B,\per}(y)} 
 \quad \text{ and } \\  
\, &T_j^{\bY_y,B_y} \, 
 = \, T_j^{B,\per,y}
\, \text{ for all }0 \le j \le m(\per)
\end{aligned} \right. \end{equation}
holds by definition. 
Furthermore we observe that for all  $y \in (\RS)^{k}$  we have 
\begin{equation} \label{trapo}  
y \in \akb \quad \Leftrightarrow \quad T^{B,\per}(y) \in \takb.
\end{equation}
Here ``$\Rightarrow $'' holds by \eqref{at} and \eqref{tildegleich} 
from the proof of Lemma~\ref{leinv}. For  ``$\Leftarrow $''
let  $y \in  (\RS)^{k}$ such that $T^{B,\per}(y) \in \takb$ and
let $(Y',B'):= \tTn(\bY, B)_{T^{B,\per}(y)}$, 
where $\tTn$  is the inverse of  $\Tn$ as defined in the last subsection. 
By induction on $j$ we can show 
\[
\alle 0 \le j \le m(\per): \quad   T^{Y',B'}_j =  T_{j}^{B,\per,y},  \;  
\per_{j,y}  = P^{Y',B'}_j, \; 
\text{ and } \; \per^{B}_{j,y} = C^{Y',B'}_j.
\]
In the inductive step $j-1 \to j$ the first assertion follows from the
induction hypothesis, the third follows from the second, and the second
follows from the first, the bijectivity of $T^{Y',B'}_j$, and 
$T_{j}^{B,\per,y}(\per_{j,y}) = T^{Y',B'}_j(P^{Y',B'}_j)$, which holds as 
\[
T_{j}^{B,\per,y}(\per_{j,y}) \, 
=  \, \per_{j,T^{B,\per}(y)} \,  
= \, \tilde{P}_j^{(\bY, B)_{T^{B,\per}(y)}}\,  
= \, P_j^{Y',B'} + \tau_j^{Y',B'} \, = \, T^{Y',B'}_j(P^{Y',B'}_j),
\]
where we have used the definition of $\takb$ 
and \eqref{tildegleich} from the proof of Lemma~\ref{leinv}. 
This completes the proof of the above assertion and we conclude that 
$(\bY_y,B_y)= (Y',B')$, which implies 
$\per_{j,y}  = P^{Y',B'}_j = P_j^{\bY_y,B_y}$ and thus \eqref{trapo}. 
Defining $g(y) :=  1_{\takb}(y) f(\bY_y,B_y)$, 
\eqref{at} and \eqref{trapo} imply 
\[
I(k,B,\per) \, = \, \Big[ \prod_{j=0}^{m(\per)} 
   \Big( \prod_{i \in \per^B_j\cap \N_k} \int dy_i \Big) 
  \Big( \prod_{i' \in \per_{j} \cap \N_k} 
 \big|1 + \partial_{\einh} t_j^{B,\per,y}(y_{i'}) \big|\Big) \Big] \,   g(T^{B,\per}(y)),  
\]
where we have also inserted the definition  of $\ph$ \eqref{densdef}.
Now we transform the integrals. For $j=m(\per)$ to $0$ and 
$i \in \per^B_j \cap \N_k$ we substitute $y_i' := T_{j(i)}^{B,\per,y} y_i$. 
For $i \in \per^B_j \weg \per_j$ $~T_{j(i)}^{B,\per,y}$ is a translation by a 
constant vector, so $dy_i' = dy_i$. 
For $i \in \per_j$ the transformation only concerns the first spatial 
component of $y_i$, and Lemma~\ref{lelip} implies 
$dy_{i}' \, = \, \big| 1 + \partial_{\einh} t_j^{B,\per,y} (y_{i})\big| dy_{i}$. So 
\[
I(k,B,\per) \, = \,  \Big[\prod_{j=0}^{m(\per)} 
    \Big(\prod_{i \in \per^B_j \cap \N_k} \int dy'_i \Big)\Big] \, g(y') \, 
=  \, \int dy \, 1_{\takb}(y) f(\bY_y,B_y),
\]
and we are done as the same arguments show that the right hand side of 
\eqref{dens} equals 
\[ 
 e^{-4 n^2} \sum_{k \ge 0} \frac{1}{k!}
\sideset{}{'}\sum_{B \subset E_n(\bY^k)} \sum_{\per \in \Pi_k} 
 \int dy \,   1_{\takb}(y) f(\bY_y,B_y).
\]
Analogously the density function can be shown to be well defined: 
For all $\bY \in \Y$  
\[
\begin{split}
&\nu_{\Lan} \otimes \pinh(``\ph \text{ well defined}''|\bY)
= \,  e^{-4 n^2} \sum_{k \ge 0} \frac{1}{k!} \sideset{}{'}
 \sum_{B \subset E_n(\bY^k)}\sum_{\per \in \Pi_k} I'(k,B,\per), \text{ with}\\
&I'(k,B,\per) :=  \Big[\prod_{j=0}^{m(\per)}  
 \Big( \prod_{i \in \per^B_j \cap \N_k} \int dy_i \Big) \Big( \prod_{i' \in \per_{j}\cap \N_k}  
 1_{\{ \partial_{\einh} t_j^{B,\per,y}(y_{i'}) \text{ exists} \}} \Big) \Big]  1_{\akb}(y)  
.
\end{split}
\]
As $t_j^{B,\per,y}$ is $\einh$-hLd,
we have for arbitrary  $r \in \R$, $\si \in S$, $k$, $B$, $\per$, and $y$ 
as above that $\partial_{\einh} t_j^{B,\per,y}(.,r,\si)$ exists a.s..
So we may replace all indicator functions in the above product by $1$ using
Fubini's theorem, and the above probability equals $1$. 


\subsection{Key estimates: Lemma~\ref{ngross}}

For \eqref{dichtenklein} let $(Y,B) \in \Gn$. By Lemma~\ref{lelip} we have $|\partial_{\einh} \tyb_k(.)| \le 1/2$.
Using the inequality $-\log(1-a) \le 2a$, which holds for 
$0 \le a \le 1/2$, we thus obtain
\[ \begin{split} 
 f_{R,n}&(Y,B) \, := - \log \ph(Y,B)  -  \log \bph(Y,B) \\
&= \, - \sum_{k = 0}^{m(Y,B)}  \sum_{y \in \pyb_k}
 \log\big(1 - (\partial_{\einh} \tyb_k(y))^2\big) \, 
 \le \, \sum_{k = 0}^{m(Y,B)}  \sum_{y \in \pyb_k} 
  2 (\partial_{\einh} \tyb_k(y))^2.
\end{split} \]
By Lemma~\ref{lelip} $\partial_{\einh} \tyb_k(y)$ equals either $0$ or 
$\partial_{\einh} \tyb_0(y)$ or  $\partial_{\einh} m_{y',\tyb(y')}(y)$ 
for some $y' \in Y$ with $(y',y) \in \Kopp$. 
Using \eqref{abdeftx} we see that 
\begin{displaymath}
\|\partial_{\einh} m_{y',\tyb(y')}\|  
\le  (\tn(|y'|-\cK) -\tyb(y')) \cf, 
\end{displaymath}
\[
\text{ where } \quad \tyb(y') \, \ge \, \ayb(y') \, 
= \, \bigwedge_{y'' \in C_{Y,B_+}} \tn(|y''|).
\]
Furthermore for $y = (r_1,r_2,\si) \in \RS$ we have       
\begin{displaymath}
|\partial_{\einh} \tyb_0(y)| \, = \, 1_{\{n \ge |r_1| > |r_2| \vee R\}}
\tau \frac{q(|r_1| - R)}{Q(n-R)} \, 
\le \,  1_{\{y \in \Lan\}} \tau \frac{q(|y| - R)}{Q(n-R)},
\end{displaymath}
so we can estimate $f_{R,n}(Y,B)$ by the sum of the following two expressions: 
\begin{equation} \begin{split} \label{sig4}
\Sigma_2(R,n,Y) \,  &:= \,  2\tau^2 \sideset{}{}\sum_{y \in Y} 
  1_{\{y \in \Lan\}} \frac{q(|y| -R)^2}{Q(n-R)^2},\\
\Sigma_3(R,n,Y,B) \, &:=   \, 2 \cf^2 
  \sum_{y,y',y'' \in Y}  1_{\Kopp}(y',y)
   1_{\{y' \stackrel{Y,B_+}{\longleftrightarrow} y''\}} \Taun(y',y''),
\end{split} \end{equation}
where we have used the shorthand notation \eqref{Taun}. 
Using these terms in the definition \eqref{good} of $\Gn$ we get 
\eqref{dichtenklein}. 
For a proof of \eqref{taylor} we first note that 
for all  $y,y' \in \RS$, $\vartheta \in [-1,1]$ with   
$(y,y' + s\einh) \notin \Ko$ for all $s \in [-\vartheta,\vartheta]$ we 
can estimate $\bU(y,y'+\vartheta \einh) + \bU(y,y' - \vartheta \einh) 
  - 2 \bU(y,y')$ by
\[
 \varphi^{\bU}_{y,y'}(\vartheta) +  \varphi^{\bU}_{y,y'}(-\vartheta) 
  - 2 \varphi^{\bU}_{y,y'}(0)\, \le \, \sup_{s \in [-\vartheta,\vartheta]} 
  \frac{d^2}{dt^2} \varphi^{\bU}_{y,y'}(s) \vartheta^2 \, 
\le \, \psi(y,y') \vartheta^2,
\]
using Taylor expansion of  $\varphi^{\bU}_{y,y'}$ and 
the $\psi$-domination of the $\einh$-derivatives.
Now let $(Y,B) \in \Gn$. W.l.o.g. we may assume that the right hand side of  
\eqref{taylor} is finite. Introducing 
$\vartheta_{y,y'} :=  \tyb(y') - \tyb(y)$  for $y,y' \in \en(Y)$ we have 
\begin{displaymath}
\begin{split}
&H^{\bU}_{\La_{n}}(\Tnb Y)\,  
 + \,  H^{\bU}_{\La_{n}}(\iTnB Y) \,  
 - \, 2H^{\bU}_{\La_{n}}(Y)\\ 
& \,= \sum_{yy' \in \en(Y)} 
     [ \bU(y,y' + \vartheta_{y,y'}\einh) 
 + \bU(y,y' - \vartheta_{y,y'}\einh) - 2\bU(y,y')]\\  
& \, \le \sum_{yy' \in \en(Y)} \psi(y,y') 
 \, (\tyb(y) - \tyb(y'))^2  \quad 
=: \quad f_{R,n}(Y,B).
\end{split}
\end{displaymath}
In the first step we have used that  $\bU$ is $\einh$-invariant, and 
in the second step we are allowed to apply the above Taylor estimate 
as for  $(y,y') \notin \Ko$ we have  $(y,y' + s \einh) \notin \Ko$ 
for all $s \in [-\vartheta_{y,y'} ,\vartheta_{y,y'} ]$
by \eqref{grepss}, and for $(y,y') \in \Ko$ we have $ \vartheta_{y,y'} = 0$
by \eqref{vercon}. 
Applying the arithmetic-quadratic mean inequality to 
\[
\Big((\tyb(y) - \tn(|y|)) + (\tn(|y|) -\tn(|y'|)) 
 + (\tn(|y'|)- \tyb(y')) \Big)^2 
\]
we obtain
\begin{displaymath}
\begin{split}
f_{R,n}(Y,B) \,
&\le \, 6 \sideset{}{^{\neq}}\sum_{y,y' \in Y} 
 \psi(y,y') \, ( \tn(|y|) - \tyb (y) )^2 \\ 
&\qquad + 3 \sideset{}{^{\neq}}\sum_{y,y' \in Y}
 1_{\{|y| \le |y'|\}}\, \psi (y,y') \, (\tn(|y|)-\tn(|y'|))^2. 
\end{split}
\end{displaymath}
In the first sum on the right hand side we again use \eqref{abdeftx} to estimate
\[ 
(\tn(|y|) - \tyb (y) )^2 \, 
\le \, \sum_{y'' \in Y} 1_{\{|y| \le |y''|\}} 
1_{\{y \stackrel{Y,B_+}{\longleftrightarrow} y''\}} 
(\tn(|y|) -\tn(|y''|))^2, 
\]
so  $f_{R,n}(Y,B)$ can be estimated by 
the sum of the two following expressions:
\begin{equation}\begin{split} \label{sig2}
\Sigma_4(R,n,Y) \, &:= \, 3  \sideset{}{^{\neq}} \sum \limits_{y,y' \in Y} 
  \psi(y,y')  \Taun(y,y'),\\
\Sigma_5(R,n,Y,B) \, &:= \, 6  \sum \limits_{y,y',y'' \in Y}
  1_{\{y \stackrel{Y,B_+}{\longleftrightarrow} y''\}}
  \psi(y,y')  \Taun(y,y'').
\end{split} \end{equation}
Inserting these sums into the definition of $\Gn$ in  \eqref{good},
we obtain \eqref{taylor}.


\subsection{Set of good configurations: Lemma~\ref{lebadsmall}}

The set of good configurations $\Gn$ is defined in terms of 
the cluster range $\ryb$ and the functions  
$\Sigma_i(R,n,Y,B)$, see \eqref{good}. We will show that the 
$\mu \otimes \pin$-expectation of 
$\ryb$ is finite and independent of $n$ and that for fixed 
$R$ the expectation of every $\Sigma_i(R,n,Y,B)$ tends to $0$ for $n \to \infty$. Then Markov's inequality implies the desired result: 
We can first choose $R > n'$ such that  
$\mu \otimes \pin(\ryb \ge  R) < \de/2$ for all $n$, 
and we may then choose an $n >R$ such that 
$\mu \otimes \pin( \sum_{i=1}^5 \Sigma_i(R,n,Y,B) \ge \de) < \de /2$.\\

For $Y \in \Y$,  $B \subset E(Y)$ and any path $y_0,...,y_m$ in the graph 
$(Y,B_+)$ such that $y_0 \in \La_{n'}$ we have  
$|y_m| \, \le \, n' + \sum_{k=1}^m |y_k-y_{k-1}|$. 
By considering all possibilities for such paths we thus obtain
\[
\ryb \, \le \, n' +  \sum_{m \ge 1} 
 \quad \sideset{}{^{\neq}} \sum_{y_0,\ldots ,y_m \in Y}   1_{\{y_0 \in \La_{n'} \}}  
  \prod_{i=1}^m   1_{ \{y_i y_{i-1} \in B_+\} } \sum_{k = 1}^m |y_k-y_{k-1}|.
\]
Under the  Bernoulli measure  $\pin(dB|Y)$ the events 
$\{y_i y_{i-1} \in B_+ \}$ are independent, and for $g: (\RS)^2 \to \R$, 
$g := 1_{\Kopp \weg \KU} +  \tiu$ we have 
\begin{equation} \label{pig}
\int \pin(dB|Y)  1_{ \{y_i y_{i-1} \in B_+\} }  \, 
\le \,1_{\KU}(y_{i-1},y_i) +g(y_{i-1},y_i). 
\end{equation} 
Using the hard core property  \eqref{hardcore} and  Lemma~\ref{lekorab}
we obtain
\[ \begin{split}
r \, &:= \, \int \mu(dY) \int \pin(dB|Y) \ryb \, - \, n'\\
&\le \, \sum_{m \ge 1} \sum_{k = 1}^m \int \mu(dY) \, \sideset{}{^{\neq}}
 \sum_{y_0,\ldots ,y_m \in Y} 1_{\{y_0 \in \La_{n'} \}}
   |y_k-y_{k-1}| \prod_{i=1}^m g(y_{i-1},y_i)\\
&\le \, \sum_{m \ge 1} \sum_{k = 1}^m \, (z\xi)^{m+1} \, 
   \int dy_0 \ldots  dy_{m}   1_{\{y_0 \in \La_{n'} \}}
   |y_k-y_{k-1}| \prod_{i=1}^m g(y_{i-1},y_i).
\end{split} \] 
Setting $c_g := (1 + \cK^2)\ex  + \cu$, we conclude from \eqref{dec} 
and \eqref{bpsi} that we have 
\begin{equation} \label{intg}
\begin{split}
&\int g(y,y') \, dy' \, \le \, \ex  \quad \text{ and }\\
& \int g(y,y')|y-y'|\, dy' \, 
\le \,  \int g(y,y')(1 + |y-y'|^2)\, dy' \, 
\le \, c_g
\end{split} \end{equation}
for all $y \in \RS$,  
hence we can estimate the integrals over $dy_i$ in the above expression 
beginning with  $i =m$. These estimates give $m-1$ times a factor 
$\ex$ and one time a factor $c_g$. Finally the integration over  $dy_0$ 
gives an additional factor  $\la^2(\La_{n'}) = (2n')^2$. Thus   
\[
r\, \le  \, (2n'z\xi)^2 c_g \sum_{m \ge 1}  
  m (\ex z\xi )^{m-1} \, < \, \infty, \quad \text{ as } \ex z\xi < 1.
\]
This gives the finiteness of the expectation of the cluster range. 
%
%
%
%
The functions $\Sigma_i(R,n,Y,B)$  have been specified in 
\eqref{sig1}, \eqref{sig4} and \eqref{sig2}:
\[ \begin{split}
& \Sigma_1  =   4 \cf^2 \sideset{}{^{}}\sum \limits_{y,y' \in Y} 
  1_{\{y \stackrel{Y,B_+}{\longleftrightarrow} y'\}}
  \Taun(y,y'), \quad  
 \Sigma_4 = 3  \sideset{}{^{\neq}} \sum \limits_{y,y' \in Y} 
  \psi(y,y') \Taun(y,y'),  \\
& \Sigma_2  =  2\tau^2 \sideset{}{}\sum \limits_{y \in Y}  1_{\{y \in \Lan\}}
 \frac{q(|y| -R)^2}{Q(n-R)^2}, \quad  
 \Sigma_5  =  6  \!\!\! \sum \limits_{y,y',y'' \in Y} \!\!\!
  1_{\{y \stackrel{Y,B_+}{\longleftrightarrow} y'\}} 
   \psi(y,y'') \Taun(y,y'), \\
& \Sigma_3  =  2 \cf^2  \!\!\! \sum \limits_{y,y',y'' \in Y}   \!\!\! 
 1_{\{y \stackrel{Y,B_+}{\longleftrightarrow} y'\}} 1_{\Kopp}(y,y'')
 \Taun(y,y'). 
\end{split} \]
We start with an estimate  on $\Taun$. For $s' > s$ such that 
$s' > R$ and $s < n$,
\begin{displaymath}
0 \,\le \, r(s-R,n-R) - r(s'-R,n-R) \, 
= \,  \int_{R \vee s}^{s' \wedge n}
 \frac{q(t-R )}{Q(n-R)} dt \, 
\le \,  (s'-s) \, \frac{q(s - R)}{Q(n-R)}
\end{displaymath}
by the monotonicity of  $q$. Defining $\bn := n + \cK$ and $\bR := R + \cK$ 
we thus have 
\begin{equation} 
\label{taugegenq}
\Taun(y,y') \,
\le \, 1_{\{y \in \La_{\bn} \}} \tau^2 \, (|y'|-|y|+ \cK)^2 \, 
    \frac{q(|y| - \bR)^2}{Q(\bn-\bR)^2} \quad \text{ for } y,y' \in \RS, 
\end{equation} 
using the substitution $s' := |y'|$ and $s := |y| - \cK$. (If  
$s' \le R$ or $s \ge n$ then $\Taun(y,y') = 0$.)
The following relations will give us control over the relevant terms 
of the right hand side of \eqref{taugegenq}. For $\bn \ge 2 \bR$ 
we have 
\begin{displaymath}
\begin{split}
\int_{\La_{\bn}} &dy  \, q(|y| - \bR)^2 \quad 
\le \quad \int_0^{2\bR} ds \, 8s \, 
  + \, \int_{\bR}^{\bn-\bR} ds \, 8(s+\bR) q(s)^2 \\
&\le  \, 16\bR^2 \,+ \, 32 \int_0^{\bn-\bR} q(s)ds  \quad 
\le \quad 16\bR^2 \,+ \, 32 Q(\bn-\bR).
\end{split}
\end{displaymath}
In the first step we used $q \le 1$, and in the second step
$\bR \le s$ and $sq(s) \le 2$. As 
$\lim \limits_{n \to \infty}Q(n) = \infty$ by 
$\log \log n \le Q(n)$ for $n >1$, the above implies 
\begin{equation} \label{ntoinf}
\lim_{n \to \infty} c(R,n) \, = \, 0 \quad \text{ for } \quad 
c(R,n) \, := \, \int_{\La_{\bn}} dy  \, 
\frac{q(|y| -  \bR)^2}{Q(\bn-\bR)^2}.
\end{equation} 
Finally, for $y_0,\ldots ,y_m \in \RS$ we deduce from the 
triangle inequality that 
\[
\Big| |y_m|-|y_0| + \cK \Big| \, 
\le \,  m \bigvee_{k=1}^m  |y_k-y_{k-1}| + \cK  \, 
\le \,  (m+1)( 1 \vee \cK)\Big( 1 \vee \bigvee_{k =1}^m |y_k-y_{k-1}| \Big), 
\]
\begin{equation} \label{sqmax}
\text{ so } \quad (|y_m|-|y_0| + \cK)^2 \, 
    \le \,  (m+1)^2( 1 \vee \cK^2 ) \bigvee_{k =1}^m (1 \vee |y_k-y_{k-1}|^2 ). 
\end{equation}
No we will proceed as in the first part of this section: 
For $Y \in \Y$ and $B \subset E(Y)$ we can estimate the summands of 
$\Sigma_1(R,n,Y,B)$ by considering all paths $y_0,\ldots ,y_m$ 
in the graph $(Y,B_{+})$ connecting $y=y_0$ and $y'=y_m$. By \eqref{taugegenq} 
and \eqref{sqmax} we can estimate $\Sigma_1(R,n,Y,B)$ by a constant $c$ times 
\[
\sum_{m \ge 0} (m+1)^2 \sum_{k=1}^m \quad \sideset{}{^{\neq}}
 \sum_{y_0,\ldots ,y_m \in Y}   1_{\{y_0 \in \La_{\bn} \}}
 \frac{q(|y_0|-\bR)^2}{Q(\bn-\bR)^2} ( 1 \vee |y_k-y_{k-1}|^2) 
\prod_{i=1}^m 1_{\{y_i y_{i-1} \in B_+\}}.
\]
The expectation of the last term can be estimated using \eqref{pig}, 
Lemma~\ref{lekorab} and \eqref{pig}, and thus we get the following 
upper bound for the expectation of $\Sigma_1$: 
\[ 
 (z\xi)^2 c_g c \sum_{m \ge 0}  (m+1)^3 (\ex z \xi )^{m-1} c(R,n).
\] 
Similarly we estimate the summands of $\Sigma_5(R,n,Y,B)$ by 
considering all paths $y_0,\ldots ,y_m$ in the graph 
$(Y,B_{+})$ connecting $y=y_0$ and $y'=y_m$ and by distinguishing the 
cases $y_j = y''$ and $y_j \neq y'' \alle j$. By \eqref{taugegenq} and 
\eqref{sqmax} we can estimate $\Sigma_5(R,n,Y,B)$ by a constant $c$ times  
\[ 
\begin{split}
&\sum_{m \ge 0} (m+1)^2 \sum_{k = 1}^m \quad \sideset{}{^{\neq}}
 \sum_{y_0,\ldots ,y_m \in Y}   1_{\{y_0 \in \La_{\bn} \}}
 \frac{q(|y_0|-\bR)^2}{Q(\bn-\bR)^2} ( 1 \vee |y_k-y_{k-1}|^2) \\
& \qquad \times \prod_{i=1}^m 1_{\{y_i y_{i-1} \in B_+\}} 
 \Big[ \sum_{y'' \in Y , y'' \neq y_j \alle j} \psi(y_0,y'')
 +  \sum_{j=0}^{m}  \psi(y_0,y_j)\Big].  
\end{split} \]
We estimate the second sum in the brackets by  $\cpsi (m+1)$ using \eqref{bpsi}. 
Proceeding as above we see that the expectation of $\Sigma_5(R,n,Y,B)$ 
can be estimated by 
\[
(z\xi)^2 c_g c \sum_{m \ge 0} (m+1)^3  (\ex z\xi )^{m-1} 
 c(R,n) \Big(  z\xi \cpsi  +  \cpsi (m+1)  \Big). 
\]
We estimate $\Sigma_3(R,n,Y,B)$ by a similar term as the one  for $\Sigma_5$, where the function $\psi$ is now replaced by $1_{\Kopp}$, so 
the expectation of $\Sigma_3(R,n,Y,B)$ can be estimated by 
\[
 (z \xi)^2 c_g c \sum_{m \ge 0}  (m+1)^3 (\ex z \xi )^{m-1} \,   
 c(R,n)(z \xi \ex + m+1). 
\]
Analogously we can estimate the expectation of $\Sigma_4(R,n,Y,B)$ by 
\[ 
 c (z \xi)^2 \int_{\La_{\bn}} dy \, \frac{q(|y| -  \bR)^2}{Q(\bn-\bR)^2}
     \int dy'  \psi(y,y')(1 \vee |y-y'|^2) \,  
\le \,   c (z \xi)^2 c(R,n) \cpsi
\]
for some constant $c$, and finally the expectation of $\Sigma_2(R,n,Y)$ can be estimated by $2 z \xi  \tau^2 c(R,n)$. 
In the bounds of the expectations of the above terms the sums over $m$ 
have finite values by \eqref{dec}, so we are done by \eqref{ntoinf}. 
\end{sloppypar}

\bigskip
\bigskip

\noindent
{\bf Acknowledgements:}\\

\noindent
I would like to thank H.-O. Georgii for suggesting the problem 
and many helpful discussions and F.~Merkl for helpful comments.


\renewcommand{\thesection}{}

\setlength{\parindent}{0cm}

\end{document}